\documentclass[11pt]{article}
\usepackage{dsfont} 
\usepackage{amssymb} 
\usepackage{graphicx}
\usepackage{extarrows}

\usepackage{mathrsfs} 

\begin{document}
\title{{\bf\Large Approximating dynamics of a singularly perturbed stochastic wave equation
with a random dynamical boundary condition
\thanks{This work was supported by the National Science
Foundation of China (Grants Nos. 10901115, 11071177, 10971225 and 11028102), the NSF Grant  1025422,      the
Scientific Research Fund of Science and Technology Bureau of
Sichuan Province (Grant Nos. 2012JQ0041 and 2010JY0057), and the Fundamental Research Funds for the Central Universities (HUST  No. 2010ZD037).  } } }

\author{Guanggan Chen\\
College of Mathematics and Software Science,\\
Sichuan Normal University,
Chengdu, 610068, China\\
\emph{E-mail: chenguanggan@hotmail.com}\\
\\
Jinqiao Duan\\
Institute for Pure and Applied Mathematics,\\
University of California,
Los Angeles, CA 90095, USA\\
\emph{E-mail: jduan@ipam.ucla.edu }\\ \& \\
Department of Applied Mathematics, \\
Illinois Institute of Technology,
Chicago, IL 60616, USA \\
\emph{E-mail: duan@iit.edu }\\
\\
Jian Zhang\\
College of Mathematics and Software Science, \\
Sichuan Normal University,
Chengdu, 610068, China\\
\emph{E-mail: zhangjiancdv@sina.com}
 }

\date{\today}
\maketitle
\baselineskip 6.8mm
\begin{center}
\begin{minipage}{120mm}
{\bf \large Abstract:} {\small This work is concerned with a
singularly perturbed stochastic nonlinear wave equation with a
random dynamical boundary condition. A splitting skill is used to
derive the approximating equation of the system in the sense of
probability distribution, when the singular perturbation parameter is
sufficiently small. The approximating equation is a stochastic
parabolic equation
  when the power exponent of singular
perturbation parameter is in $[1/2, 1)$, but   a deterministic hyperbolic (wave)
equation     when  the power exponent   is in $(1, +\infty)$. }
\end{minipage}
\end{center}

\renewcommand{\theequation}{\thesection.\arabic{equation}}
\setcounter{equation}{0}

\vspace{0.5cm}

\section{ Introduction }

\quad\quad Stochastic nonlinear wave equations play an important
role in describing the propagation of waves in certain systems or
media,  such as atmosphere, oceans, sonic booms, traffic flows,
optic devices and quantum fields,  when random fluctuations are
taken into account (\cite{Chen1, Chow1, DZ1, Pardoux, RS, Whitham}).
They have been studied   recently by  a number of authors(see \cite{Chen2, Chen3, Chow2,
FW, LS07, M, ZYO}). For some wave systems
on bounded domains, noise may affect the system evolution through
the boundary in terms of random boundary conditions. Dirichlet, Neumann and Robin boundary conditions are
{\it static boundary conditions}, because they are not involved with
time derivatives of the system state variables. On the contrary,
{\it dynamical boundary conditions} contain time derivatives of the
system state variables and arise in many physical problems
(\cite{FG, GS, PR}).
\par
In this paper, we investigate a singularly perturbed stochastic
nonlinear wave equation with a random dynamical boundary condition
\begin{equation}\label{Eq1}
\left\{
\begin{array}{ll}
\varepsilon u^\varepsilon_{tt}+u^\varepsilon_t-\bigtriangleup
u^\varepsilon+u^\varepsilon-f(u^\varepsilon)=\varepsilon^\alpha\dot{W_1}
&\quad \hbox{in}\; D,\\
\varepsilon\delta^\varepsilon_{tt}+\delta^\varepsilon_t+\delta^\varepsilon
=-u^\varepsilon_t+\varepsilon^\alpha\dot{W_2} &\quad \hbox{on}\;
\partial D,\\
\delta^\varepsilon_t=\frac{\partial u^\varepsilon}{\partial {\bf n}}
&\quad \hbox{on}\; \partial D,\\
u^\varepsilon(0)=u_0, u^\varepsilon_t(0)=u_1,
\delta^\varepsilon(0)=\delta_0, \delta^\varepsilon_t(0)=\delta_1.
\end{array}
\right.
\end{equation}
Here $u^\varepsilon(x, t)$ is the unknown wave amplitude,
$\varepsilon $ is a small positive singular perturbation parameter
($0 <\varepsilon \ll 1$), and the power exponent $\alpha$ is in
$[1/2, 1)$ or $ (1,+\infty)$. Moreover, $W_1$ and $W_2$ are two
independent Wiener processes, which will be specified in details in
the next section. The symbol $\frac{\partial }{\partial {\bf n}}$
denotes the unit outer normal derivative on the boundary $\partial
D$ of a bounded domain $D$ in $\mathbb{R}^3$. Note that
$\delta^\varepsilon_t$ is the outer normal derivative of
$u^\varepsilon$ on the boundary. We often write $u^\varepsilon(x,
t)$ as $u^\varepsilon(t)$. In particular, in this paper we will only
concern with the case of the nonlinear term $f(u^\varepsilon)=\sin
u^\varepsilon$ (the Sine-Gordon equation).
\par
The system (\ref{Eq1}) arises in the modeling of gas dynamics in an
open bounded domain $D$, with points on boundary acting like a
spring reacting to the excess pressure of the gas (see \cite{MI}).
  Chen and Zhang \cite{Chen2} studied the
long time behavior of the solutions of the system (\ref{Eq1})
without the singular perturbation parameter. Also Chen, Duan and
Zhang \cite{Chen3} derived the effective dynamics of the system
(\ref{Eq1}) on a bounded domain perforated with small holes. For the deterministic case of the system (\ref{Eq1}), Beale
\cite{Beale1, Beale2} and Mugnolo \cite{Mugnolo} established the
well-posedness   in some
special cases. Cousin, Frota and Larkin \cite{CFL} studied the
global solvability and asymptotic behavior. Frigeri \cite{Frigeri}
considered   large time dynamical behavior.
\par
The singular perturbation issues of wave equations  have been
studied extensively. On the one hand, for  a deterministic wave
equation with a static boundary condition, Hale and Raugel \cite{Hale-R}
and Mora \cite{Mora} studied the approximation as the perturbation
parameter $\varepsilon$ goes to zero.   For   deterministic wave
equations with dynamical boundary conditions, Rodriguez-Bernal and
Zuazua \cite{RZ-1, RZ-2} and Popsescu and Rodriguez-Bernal \cite{PR}
considered the singular limiting equations. On the other hand, for   stochastic
wave equations with homogeneous boundary conditions, Cerrai and
Freidlin \cite{CF1, CF2}, Lv, Roberts and Wang \cite{LR, LW, WL, WLR}
investigated the approximation as the perturbation parameter
$\varepsilon$ goes to zero.

\par
In the present paper, we investigate the singular
perturbations of the stochastic wave equation with   random dynamical
boundary conditions. Our goal is to derive the approximating equation
of the system (\ref{Eq1}) for   sufficiently small
  parameter $\varepsilon$. There are two key points to
achieve this goal: The first is to establish the tightness of the
solutions, and the second is to construct
the approximating equation of the system (\ref{Eq1}).
\par
For the first key point, the tightness of solutions for the system
(\ref{Eq1}) heavily depends on the almost sure boundedness of the solutions, independent of the parameter $\varepsilon$. However,
since the parameter $\varepsilon$ disturbs the system (\ref{Eq1}),
it is difficult to derive the almost sure boundedness independent of
    $\varepsilon$. As showing in Chen and Zhang
\cite{Chen2}, the classic energy relation of this stochastic system
(\ref{Eq1}) does not directly imply the a priori estimate of the
solutions. Meanwhile, as we will see, the pseudo energy argument
especially proposed in Chow \cite{Chow2} and Chen and Zhang
\cite{Chen2} for stochastic wave equations also does not lead to the a
priori estimate of the solutions. Therefore,
for the system (\ref{Eq1}), we will explore a new way to establish
the a priori estimate of solutions. By applying the a priori
estimate, we could then obtain the global well-posedness and the almost sure
boundedness independent of   $\varepsilon$,
which further implies the tightness of the solutions.
\par
For the second key point, we use a splitting skill to construct the
approximating equation. Firstly, we split the solution of the system
(\ref{Eq1}) into three parts: the solution of a linear random
ordinary differential equation (RODE), the solution of a random
partial differential equation (RPDE), and the solution of a linear
stochastic ordinary differential equation (SODE). Then we analyze
their respective approximations     for sufficiently
small $\varepsilon$. Finally, we derive the approximating equation
of the system (\ref{Eq1}) for the sufficiently small $\varepsilon$
in the sense of probability distribution, which is a stochastic
parabolic equation with a dynamical boundary condition for $\alpha\in
[1/2, 1)$,  and   a deterministic wave equation with a dynamical
boundary condition for $\alpha\in (1, +\infty)$.
\par
We especially remark that the power exponent (of the singular
perturbation parameter), $\alpha$, is in the set $[1/2, 1)\cup (1,
+\infty)$.  The case of $0<\alpha<1/2$ is not covered in our results, as the
condition of $\alpha\geq 1/2$ plays two crucial roles in our
work: One is in deriving the a priori estimate of solutions for
the system (\ref{Eq1}) (see Proposition 3.2), and the other is in
deriving the almost sure boundedness of solutions for the split
linear stochastic ordinary differential equation (see the proof of
Theorem 4.1). In addition, for the case of $\alpha=1$, in analyzing
the approximation of the decomposition of the system ({\ref{Eq1}})
(see the proof of Theorem 4.1), there is no difference of
convergence velocity between $O(\varepsilon)$ and
$O(\varepsilon^\alpha)$ as $\varepsilon$ tends to zero, which means
that the final approximating equation of the system (\ref{Eq1}) is
just itself.

\par
This paper is organized as follows. In the next section, we
present   preliminary results including the local well-posedness of the
system (\ref{Eq1}) (Proposition 2.1). In section 3, we first derive
the pseudo energy relation of the system (\ref{Eq1}) (Proposition
3.1), which   implies certain estimates, independent of
$\varepsilon$, for one part of the solution of the system (\ref{Eq1})
(Remark 3.2). With the help of  these estimates, we
  establish the estimates, independent
of   $\varepsilon$, for the other parts of the solution of the
system (\ref{Eq1}) (Proposition 3.4). All the a priori estimates
   for the solution  play an important
role in proving the global well-posedness (Proposition 3.5), the almost
sure boundedness (Remark 3.3) and the tightness (Proposition 3.8).
In section 4, we examine   the solution  as decomposed into three parts    (Proposition 4.1), and
further     derive the approximating
equation of the system (\ref{Eq1}), in the sense of probability
distribution (Theorem 4.1).
\par

\renewcommand{\theequation}{\thesection.\arabic{equation}}
\setcounter{equation}{0}

\section{Preliminaries}

\quad \quad Consider the Wiener processes $W_1(t)$ and $W_2(t)$, defined
on a complete probability space $(\Omega, \mathcal{F}, \mathbb{P})$
with a filtration $\{\mathcal{F}_t\}_{t\in \mathbb{R}}$,  and
two-sided in time with values in $L^2(D)$ and $L^2(\partial D)$,
respectively. Further  assume that $W_1(t)$ and $W_2(t)$ are
 independent and that their covariance operators,   $Q_1$ and
$Q_2$, are   symmetric nonnegative operators,  satisfying
$Tr Q_1<+\infty$ and $Tr Q_2<+\infty$. Their expansions are given as
follows
$$
\begin{array}{l}
 W_1(t)=\sum\limits_{i=1}^{+\infty}\sqrt{\alpha_{1i}}\beta_{1i} e_i,
\quad \hbox{with}\quad Q_1e_i=\alpha_{1i} e_i,\\
 W_2(t)=\sum\limits_{i=1}^{+\infty}\sqrt{\alpha_{2i}}\beta_{2i} \gamma(e_i),
\quad \hbox{with}\quad Q_2e_i=\alpha_{2i} \gamma(e_i),
\end{array}
$$
where $\{e_{i}\}_{i\in \mathbb{N}}$ is an orthonormal basis of
$L^2(D)$, and $\gamma$ is the trace operator from $D$ to $\partial
D$. Moreover, $\{\beta_{1i}\}_{i\in \mathbb{N}}$ and
$\{\beta_{2i}\}_{i\in \mathbb{N}}$ are two sequences of mutually
independent (two-sided in time) standard scalar Wiener processes in
the probability space $(\Omega, \mathcal{F}, \mathbb{P})$.
\par
Then the system (\ref{Eq1}) can be written in the It$\hat{o}$ form
as follows
\begin{equation}\label{Eq2}
\left\{
\begin{array}{ll}
du^\varepsilon=v^\varepsilon dt &\quad \hbox{in}\; D,\\
dv^\varepsilon=(\frac{1}{\varepsilon}\bigtriangleup
u^\varepsilon-\frac{1}{\varepsilon}u^\varepsilon-\frac{1}{\varepsilon}v^\varepsilon+\frac{1}{\varepsilon}\sin
u^\varepsilon)dt+\varepsilon^{\alpha-1}dW_1 &\quad \hbox{in}\;
D,\\
d\delta^\varepsilon=\theta^\varepsilon dt & \quad \hbox{on}\; \partial D, \\
d\theta^\varepsilon=(-\frac{1}{\varepsilon}\theta^\varepsilon-\frac{1}{\varepsilon}\delta^\varepsilon
-\frac{1}{\varepsilon}v^\varepsilon)dt +\varepsilon^{\alpha-1}dW_2 &
\quad \hbox{on}\;
\partial D,\\
\delta^\varepsilon_t=\frac{\partial u^\varepsilon}{\partial {\bf
n}}& \quad \hbox{on}\;
\partial D,\\
u^\varepsilon(0)=u_0, v^\varepsilon(0)=v_0=u_1,
\delta^\varepsilon(0)=\delta_0,
\theta^\varepsilon(0)=\theta_0=\delta_1.
\end{array}
\right.
\end{equation}
\par
Now we define
$$
A^\varepsilon=\left(
\begin{array}{cccc}
0&\; I&\; 0&\; 0\\
\frac{1}{\varepsilon}(\bigtriangleup-I)&\; -\frac{1}{\varepsilon}I&\; 0&\; 0\\
0&\; 0&\; 0&\; I\\
0&\; -\frac{1}{\varepsilon}I&\; -\frac{1}{\varepsilon}I&\;
-\frac{1}{\varepsilon}I
\end{array}
\right) ,  F^\varepsilon(U^\varepsilon)=\left(
\begin{array}{c}
0\\
\frac{1}{\varepsilon}\sin u^\varepsilon\\
0\\
0
\end{array}
\right), W=\left(
\begin{array}{c}
0\\
W_1\\
0\\
W_2
\end{array}
\right).
$$
Let $U^\varepsilon:=(u^\varepsilon, v^\varepsilon,
\delta^\varepsilon, \theta^\varepsilon)^T$ be in the Hilbert space
$$
\mathcal{H}: =\{U^\varepsilon\in H^1(D)\times L^2(D)\times
L^2(\partial D)\times L^2(\partial D)|\; \frac{\partial
u^\varepsilon}{\partial {\bf n}}=\theta^\varepsilon\},$$
with norm
$$
\|U^\varepsilon\|_{\mathcal{H}}^2=\|u^\varepsilon\|_{H^1(D)}^2
+\|v^\varepsilon\|_{L^2(D)}^2+\|\delta^\varepsilon\|_{L^2(\partial
D)}^2+\|\theta^\varepsilon\|_{L^2(\partial D)}^2.
$$
Here and hereafter, the superscript ``$T$" denotes the transpose for
a matrix.
\par
Thus the system (\ref{Eq2})  is further rewritten as
\begin{equation}\label{Eq3}
\left\{
\begin{array}{l}
dU^\varepsilon=A^\varepsilon U^\varepsilon dt+F^\varepsilon(U^\varepsilon)dt+ \varepsilon^{\alpha-1}dW(t),\\
U^\varepsilon(0)=U^\varepsilon_0=(u_0,v_0,\delta_0,\theta_0)^T.
\end{array}
\right.
\end{equation}
For the Cauchy problem (\ref{Eq3}), it follows from Frigeri
\cite{Frigeri} that the operator $A^\varepsilon$ generates a
strongly continuous semigroup $S(t)=e^{A^\varepsilon t}$ for $t\geq
0$ on $\mathcal{H}$. Then  Equation (\ref{Eq3}) can
be formulated in the mild sense
\begin{equation}\label{E-m}
U^\varepsilon (t)=S(t) U^\varepsilon(0)+\int_0^tS(t-s)F^\varepsilon
(U^\varepsilon(s))ds+\int_0^tS(t-s) \varepsilon^{\alpha-1}dW(s).
\end{equation}
\par
{\bf Proposition 2.1} {(\bf Local well-posedness)}\quad {\it Let the
initial datum $U^\varepsilon(0)$ be a $ \mathcal{F}_0$-measurable
random variable with value in $\mathcal{H} $. Then the Cauchy
problem (\ref{Eq3}) has a unique local mild solution
$U^\varepsilon(t)$ in $C([0, \tau^*), \mathcal{H} )$, where $\tau^*$
is a stopping time depending on $U^\varepsilon(0)$ and $\omega$.
Moreover, the mild solution $U^\varepsilon(t)$ is also a weak
solution in the following sense
\begin{equation}
\langle U^\varepsilon (t), \phi\rangle_{\mathcal{H} }= \langle
U^\varepsilon(0), \phi \rangle_{\mathcal{H} } +\int_0^t \langle
A^\varepsilon U^\varepsilon (s), \phi\rangle_{\mathcal{H} }ds+
\int_0^t \langle F^\varepsilon(U^\varepsilon (s)),
\phi\rangle_{\mathcal{H} }ds+\int_0^t \langle
\varepsilon^{\alpha-1}dW(s),\phi \rangle_{\mathcal{H} }
\end{equation}
for any $t\in[0, \tau^*)$ and $\phi\in \mathcal{H} $. }
\par
Using the cut-off function method and combining with Theorem 7.4  and the stochastic Fubini theorem
in \cite{DZ1}, we  can prove Proposition 2.1. Please refer to Chen
and Zhang\cite{Chen2}.

\par

\renewcommand{\theequation}{\thesection.\arabic{equation}}
\setcounter{equation}{0}

\section{Boundedness and tightness}

\quad\quad In this section, we will establish the almost sure
boundedness independent of the parameter $\varepsilon$ and the
tightness of solutions for the system (\ref{Eq1}). Due to the singular perturbation in the system (\ref{Eq1}), the classic
energy method and the pseudo energy method does not directly imply
the almost sure boundedness, independent of the parameter
$\varepsilon$, of solutions. We will explore a new way to do it.
\par
For  a real parameter $r$ in $(0,1)$, we define
\begin{equation}\label{T-Psolution}
v_r^\varepsilon=v^\varepsilon+ru^\varepsilon\quad \hbox{and}\quad
\theta_r^\varepsilon=\theta^\varepsilon+r\delta^\varepsilon,
\end{equation}
with $(u^\varepsilon, v^\varepsilon, \delta^\varepsilon,
\theta^\varepsilon)^T$ being the solution of the Cauchy problem
(\ref{Eq2}). Then     $U_r^\varepsilon= (u^\varepsilon,
v_r^\varepsilon, \delta^\varepsilon, \theta_r^\varepsilon)^T \in
\mathcal{H}$ satisfies the following equation
\begin{equation}\label{Eq-P}
\left\{
\begin{array}{ll}
du^\varepsilon=(v_r^\varepsilon-r u^\varepsilon)dt & \hbox{in}\;
D\times
[0, \tau^*),\\
dv_r^\varepsilon=(\frac{1}{\varepsilon}\bigtriangleup
u^\varepsilon-\frac{1}{\varepsilon}(1-r+\varepsilon
r^2)u^\varepsilon-\frac{1}{\varepsilon}(1-\varepsilon
r)v_r^\varepsilon+\frac{1}{\varepsilon}\sin
u^\varepsilon)dt+\varepsilon^{\alpha-1}dW_1(t) & \hbox{in}\; D
\times
[0, \tau^*),\\
d\delta^\varepsilon=(\theta_r^\varepsilon-r\delta^\varepsilon)dt &  \hbox{on}\;
\partial D \times [0, \tau^*), \\
d\theta_r^\varepsilon=(-\frac{1}{\varepsilon}(1-\varepsilon
r)\theta_r^\varepsilon -\frac{1}{\varepsilon}(1-r+\varepsilon
r^2)\delta^\varepsilon-\frac{1}{\varepsilon}v_r^\varepsilon+\frac{r}{\varepsilon}
u^\varepsilon)dt+\varepsilon^{\alpha-1}dW_2(t) &  \hbox{on}\; \partial D\times [0, \tau^*),\\
\delta^\varepsilon_t=\frac{\partial u^\varepsilon}{\partial {\bf
n}}& \hbox{on}\;
\partial D \times [0, \tau^*),\\
u^\varepsilon(0)=u_0, v_r^\varepsilon(0)=v_0+r u_0:=v_{r0},\\
\delta^\varepsilon(0)=\delta_0,
\theta_r^\varepsilon(0)=\theta_0+r\delta_0:=\theta_{r0}.
\end{array}
\right.
\end{equation}
\par

Define a pseudo energy functional $\mathcal{E}_r^\varepsilon(t)$
for the Cauchy problem (\ref{Eq3}) as follows
$$
\begin{array}{lll}
\mathcal{E}_r^\varepsilon(t)&:=&\varepsilon\|v_r^\varepsilon(t)\|_{L^2(D)}^2+\|\bigtriangledown
u^\varepsilon(t)\|_{L^2(D)}^2+(1-r+\varepsilon r^2)\|u^\varepsilon(t)\|_{L^2(D)}^2\\
&& +\varepsilon \|\theta_r^\varepsilon(t)\|_{L^2(\partial D)}^2
+(1-r+\varepsilon r^2)\|\delta^\varepsilon(t)\|_{L^2(\partial D)}^2
+4\|\cos \frac{u^\varepsilon(t)}{2}\|_{L^2(D)}^2\\
&& +2r \langle
u^\varepsilon(t),\delta^\varepsilon(t)\rangle_{L^2(\partial D)}.
\end{array}
$$
\par
{\bf Proposition 3.1 (Pseudo energy equation)}\quad{\it Let the
initial datum $U_r^\varepsilon(0)$ be a $ \mathcal{F}_0$-measurable
random variable in $L^2(\Omega, \mathcal{H} )$. Then for any time
$t\in [0, \tau^*)$, we have
\begin{equation}\label{Eq-3.1-01}
\begin{array}{ll}
\mathcal{E}_r^\varepsilon(t)=&
\mathcal{E}_r^\varepsilon(0)-\int_0^t[2(1-\varepsilon
r)\|v_r^\varepsilon\|_{L^2(D)}^2+2r\|\bigtriangledown
u^\varepsilon\|_{L^2(D)}^2+2(1-r+\varepsilon
r^2)r\|u^\varepsilon\|_{L^2(D)}^2\\
&+2(1-\varepsilon r)\|\theta_r^\varepsilon\|_{L^2(\partial
D)}^2+2(1-r+\varepsilon r^2)r\|\delta^\varepsilon\|_{L^2(\partial
D)}^2
]ds\\
&+2r\int_0^t\langle u^\varepsilon,\sin
u^\varepsilon\rangle_{L^2(D)}ds +4r\int_0^t\langle u^\varepsilon,
\theta_r^\varepsilon\rangle_{L^2(\partial D)}ds- 4r^2 \int_0^t
\langle
u^\varepsilon,\delta^\varepsilon\rangle_{L^2(\partial D)}ds\\
&+\int_0^t\langle 2v_r^\varepsilon, \varepsilon^{\alpha}
dW_1(s)\rangle_{L^2(D)}+ \int_0^t\langle 2\theta_r^\varepsilon,
\varepsilon^{\alpha} dW_2(s)\rangle_{L^2(\partial
D)}\\
&+\varepsilon^{2\alpha-1} TrQ_1\cdot t+\varepsilon^{2\alpha-1}
TrQ_2\cdot t.
\end{array}
\end{equation}
Moreover,
\begin{equation}\label{Eq-3.1-02}
\begin{array}{ll}
\mathbb{E} \mathcal{E}_r^\varepsilon(t)=& \mathbb{E}
\mathcal{E}_r^\varepsilon(0)-\int_0^t[2(1-\varepsilon r)\mathbb{E}
\|v_r^\varepsilon\|_{L^2(D)}^2+2r\mathbb{E} \|\bigtriangledown
u^\varepsilon\|_{L^2(D)}^2\\
&+2(1-r+\varepsilon r^2)r\mathbb{E} \|u^\varepsilon\|_{L^2(D)}^2
+2(1-\varepsilon r)\mathbb{E} \|\theta_r^\varepsilon\|_{L^2(\partial
D)}^2\\
&+2(1-r+\varepsilon r^2)r\mathbb{E}
\|\delta^\varepsilon\|_{L^2(\partial D)}^2 ]ds+2r\int_0^t\mathbb{E}
\langle u^\varepsilon,\sin u^\varepsilon\rangle_{L^2(D)}ds\\
&+4r\int_0^t\mathbb{E} \langle u^\varepsilon,
\theta_r^\varepsilon\rangle_{L^2(\partial D)}ds- 4r^2
\int_0^t\mathbb{E}  \langle
u^\varepsilon,\delta^\varepsilon\rangle_{L^2(\partial D)}ds\\
&+\varepsilon^{2\alpha-1} TrQ_1\cdot t+\varepsilon^{2\alpha-1}
TrQ_2\cdot t.
\end{array}
\end{equation}
}

\par
{\bf  Proof.}\quad  Firstly, we examine the second equation of
(\ref{Eq-P}). Set $M(v_r^\varepsilon):=\int_{D }
|v_r^\varepsilon|^2dx$. Then from the It$\hat{o}$ formula, we deduce
that
\begin{equation}\label{Eq-3.1-1}
\begin{array}{ll}
M(v_r^\varepsilon(t))=& M(v_r^\varepsilon(0))+\int_0^t\langle
M^{\prime}(v_r^\varepsilon),
\varepsilon^{\alpha-1}dW_1(s)\rangle_{L^2(D)}\\
&+\int_0^t\frac{1}{2}Tr[M^{\prime\prime}(v_r^\varepsilon)(\varepsilon^{\alpha-1}Q_1)^{\frac{1}{2}}(\varepsilon^{\alpha-1}Q_1^{\frac{1}{2}})^*]ds\\
&+ \int_0^t\langle
M^{\prime}(v_r^\varepsilon),(\frac{1}{\varepsilon}\bigtriangleup
u^\varepsilon-\frac{1}{\varepsilon}(1-r+\varepsilon
r^2)u^\varepsilon-\frac{1}{\varepsilon}(1-\varepsilon
r)v_r^\varepsilon+\frac{1}{\varepsilon}\sin
u^\varepsilon)\rangle_{L^2(D)}ds,
\end{array}
\end{equation}
with   $M^\prime(v_r^\varepsilon)=2v_r^\varepsilon$ and
$M^{\prime\prime}(v_r^\varepsilon)=2\varphi$ for any $\varphi$ in
$L^2(D)$. After some calculations, we   conclude that
\begin{equation}\label{Eq-3.1-2}
\begin{array}{ll}
&\langle
M^{\prime}(v_r^\varepsilon),(\frac{1}{\varepsilon}\bigtriangleup
u^\varepsilon-\frac{1}{\varepsilon}(1-r+\varepsilon
r^2)u^\varepsilon-\frac{1}{\varepsilon}(1-\varepsilon
r)v_r^\varepsilon+\frac{1}{\varepsilon}\sin
u^\varepsilon)\rangle_{L^2(D)}\\
=& -\frac{d}{ds}[\frac{1}{\varepsilon}\|\bigtriangledown
u^\varepsilon\|_{L^2(D)}^2+\frac{1}{\varepsilon}(1-r+\varepsilon
r^2)\|u^\varepsilon\|_{L^2(D)}^2+\frac{4}{\varepsilon}\|\cos
\frac{u^\varepsilon}{2}\|_{L^2(D)}^2]\\
&-[\frac{2r}{\varepsilon}\|\bigtriangledown
u^\varepsilon\|_{L^2(D)}^2+\frac{2}{\varepsilon}\cdot(1-r+\varepsilon
r^2)\cdot r\|u^\varepsilon\|_{L^2(D)}^2
+\frac{2}{\varepsilon}(1-\varepsilon r)\|v_r^\varepsilon\|_{L^2(D)}^2]\\
&+\frac{2}{\varepsilon }\langle v_r^\varepsilon,\frac{\partial
u^\varepsilon}{\partial {\bf n}}\rangle_{L^2(\partial
D)}+\frac{2r}{\varepsilon }\langle u^\varepsilon,\sin
u^\varepsilon\rangle_{L^2(D)}.
\end{array}
\end{equation}
It  further follows from (\ref{Eq-3.1-1}) and (\ref{Eq-3.1-2})
that
\begin{equation}\label{Eq-3.1-3}
\begin{array}{ll}
&\|v_r^\varepsilon(t)\|_{L^2(D)}^2+\frac{1}{\varepsilon
}\|\bigtriangledown
u^\varepsilon(t)\|_{L^2(D)}^2+\frac{1}{\varepsilon
}\cdot(1-r+\varepsilon r^2)
\|u^\varepsilon(t)\|_{L^2(D)}^2+\frac{4}{\varepsilon }\|\cos
\frac{u^\varepsilon(t)}{2}\|_{L^2(D)}^2\\
=& \|v_r^\varepsilon(0)\|_{L^2(D)}^2+\frac{1}{\varepsilon
}\|\bigtriangledown
u^\varepsilon(0)\|_{L^2(D)}^2+\frac{1}{\varepsilon
}\cdot(1-r+\varepsilon
r^2)\|u^\varepsilon(0)\|_{L^2(D)}^2+\frac{4}{\varepsilon }\|\cos
\frac{u^\varepsilon(0)}{2}\|_{L^2(D)}^2\\
&-\int_0^t[\frac{2}{\varepsilon }(1-\varepsilon
r)\|v_r^\varepsilon\|_{L^2(D)}^2+\frac{2r}{\varepsilon
}\|\bigtriangledown
u^\varepsilon\|_{L^2(D)}^2+\frac{2}{\varepsilon }(1-r+\varepsilon r^2)r\|u^\varepsilon\|_{L^2(D)}^2]ds\\
&+\frac{2}{\varepsilon }\int_0^t\langle
v_r^\varepsilon,\frac{\partial u^\varepsilon}{\partial {\bf
n}}\rangle_{L^2(\partial D)}ds+\frac{2r}{\varepsilon
}\int_0^t\langle u^\varepsilon,\sin
u^\varepsilon\rangle_{L^2(D)}ds\\
&+\int_0^t\langle 2v_r^\varepsilon, \varepsilon^{\alpha-1}
dW_1(s)\rangle_{L^2(D)}+\varepsilon^{2\alpha-2} TrQ_1\cdot t.
\end{array}
\end{equation}
\par
Secondly, we examine the fourth equation of (\ref{Eq-P}). Set
$M(\theta_r^\varepsilon) :=\int_{\partial D}
|\theta_r^\varepsilon|^2dx$. Note that
\begin{equation}\label{Eq-3.1-4}
\begin{array}{lll}
M(\theta_r^\varepsilon(t))&=&
M(\theta_r^\varepsilon(0))+\int_0^t\langle
M^{\prime}(\theta_r^\varepsilon),\varepsilon^{\alpha-1}dW_2(s)\rangle_{L^2(\partial
D)}\\
&&+\int_0^t\frac{1}{2}Tr[M^{\prime\prime}(\theta_r^\varepsilon)
(\varepsilon^{\alpha-1} Q_2)^{\frac{1}{2}}(\varepsilon^{\alpha-1}Q_2^{\frac{1}{2}})^*]ds\\
&&+ \int_0^t\langle M^{\prime}(\theta_r^\varepsilon),(
-\frac{1}{\varepsilon}(1-\varepsilon r)\theta_r^\varepsilon
-\frac{1}{\varepsilon}(1-r+\varepsilon
r^2)\delta^\varepsilon-\frac{1}{\varepsilon}v_r^\varepsilon+\frac{r}{\varepsilon}
u^\varepsilon     )\rangle_{L^2(\partial D)}ds,
\end{array}
\end{equation}
with   $M^\prime(\theta_r^\varepsilon)=2\theta_r^\varepsilon$ and
$M^{\prime\prime}(\theta_r^\varepsilon)=2\phi$ for any $\phi$ in
$L^2(\partial D)$. After some calculations, we   obtain that
\begin{equation}\label{Eq-3.1-5}
\begin{array}{ll}
&\langle
M^{\prime}(\theta_r^\varepsilon),(-\frac{1}{\varepsilon}(1-\varepsilon
r)\theta_r^\varepsilon -\frac{1}{\varepsilon}(1-r+\varepsilon
r^2)\delta^\varepsilon-\frac{1}{\varepsilon}v_r^\varepsilon+\frac{r}{\varepsilon}
u^\varepsilon)\rangle_{L^2(\partial D)}\\
=&-\frac{1}{\varepsilon }(1-r+\varepsilon
r^2)\frac{d}{ds}\|\delta^\varepsilon\|_{L^2(\partial D)}^2
-\frac{2}{\varepsilon }\cdot(1-r+\varepsilon
r^2)r\|\delta^\varepsilon\|_{L^2(\partial D)}^2-\frac{2}{\varepsilon
}(1-\varepsilon r)\|\theta_r^\varepsilon\|_{L^2(\partial
D)}^2\\
&-\frac{2}{\varepsilon }\langle \frac{\partial
u^\varepsilon}{\partial {\bf
n}},v_r^\varepsilon\rangle_{L^2(\partial D)}-\frac{2r}{\varepsilon
}\langle \delta^\varepsilon,v_r^\varepsilon\rangle_{L^2(\partial D)}
+\frac{2r}{\varepsilon }\langle \theta_r^\varepsilon,
u^\varepsilon\rangle_{L^2(\partial D)}.
\end{array}
\end{equation}
Then it   follows from (\ref{Eq-3.1-4}) and (\ref{Eq-3.1-5}) that
\begin{equation}\label{Eq-3.1-6}
\begin{array}{ll}
&\|\theta_r^\varepsilon(t)\|_{L^2(\partial D)}^2
+\frac{1}{\varepsilon }(1-r+\varepsilon
r^2)\|\delta^\varepsilon(t)\|_{L^2(\partial D)}^2\\
=&\|\theta^\varepsilon(0)\|_{L^2(\partial D)}^2
+\frac{1}{\varepsilon }(1-r+\varepsilon
r^2)\|\delta^\varepsilon(0)\|_{L^2(\partial D)}^2\\
&-\int_0^t[\frac{2}{\varepsilon }(1-\varepsilon
r)\|\theta_r^\varepsilon\|_{L^2(\partial D)}^2+\frac{2}{\varepsilon
}\cdot(1-r+\varepsilon
r^2)r\|\delta^\varepsilon\|_{L^2(\partial D)}^2]ds\\
&-\frac{2}{\varepsilon }\int_0^t\langle \frac{\partial
u^\varepsilon}{\partial {\bf
n}},v_r^\varepsilon\rangle_{L^2(\partial D)}ds-\frac{2r}{\varepsilon
}\int_0^t\langle
\delta^\varepsilon,v_r^\varepsilon\rangle_{L^2(\partial D)}ds
+\frac{2r}{\varepsilon }\int_0^t\langle \theta_r^\varepsilon,
u^\varepsilon\rangle_{L^2(\partial D)}ds\\
&+\int_0^t\langle 2\theta_r^\varepsilon, \varepsilon^{\alpha-1}
dW_2(s)\rangle_{L^2(\partial D)}+\varepsilon^{2\alpha-2} TrQ_2\cdot
t.
\end{array}
\end{equation}
\par
Thus, from (\ref{Eq-3.1-3}) and (\ref{Eq-3.1-6}), we have
\begin{equation}\label{Eq-3.1-7}
\begin{array}{ll}
&\|v_r^\varepsilon(t)\|_{L^2(D)}^2+\frac{1}{\varepsilon
}\|\bigtriangledown
u^\varepsilon(t)\|_{L^2(D)}^2+\frac{1}{\varepsilon }(1-r+\varepsilon
r^2)
\|u^\varepsilon(t)\|_{L^2(D)}^2+\|\theta_r^\varepsilon(t)\|_{L^2(\partial D)}^2\\
&+\frac{1}{\varepsilon }(1-r+\varepsilon
r^2)\|\delta^\varepsilon(t)\|_{L^2(\partial
D)}^2+\frac{4}{\varepsilon }\|\cos
\frac{u^\varepsilon(t)}{2}\|_{L^2(D)}^2\\
=& \|v_r^\varepsilon(0)\|_{L^2(D)}^2+\frac{1}{\varepsilon
}\|\bigtriangledown
u^\varepsilon(0)\|_{L^2(D)}^2+\frac{1}{\varepsilon }(1-r+\varepsilon
r^2)\|u^\varepsilon(0)\|_{L^2(D)}^2
+\|\theta^\varepsilon(0)\|_{L^2(\partial
D)}^2\\
&+\frac{1}{\varepsilon }(1-r+\varepsilon
r^2)\|\delta^\varepsilon(0)\|_{L^2(\partial
D)}^2+\frac{4}{\varepsilon }\|\cos
\frac{u^\varepsilon(0)}{2}\|_{L^2(D)}^2\\
&-\int_0^t[\frac{2}{\varepsilon }(1-\varepsilon
r)\|v_r^\varepsilon\|_{L^2(D)}^2+\frac{2r}{\varepsilon
}\|\bigtriangledown u^\varepsilon\|_{L^2(D)}^2+\frac{2}{\varepsilon
}\cdot(1-r+\varepsilon
r^2)r\|u^\varepsilon\|_{L^2(D)}^2\\
&+\frac{2}{\varepsilon }(1-\varepsilon
r)\|\theta_r^\varepsilon\|_{L^2(\partial D)}^2+\frac{2}{\varepsilon
}\cdot(1-r+\varepsilon r^2)r\|\delta^\varepsilon\|_{L^2(\partial
D)}^2
]ds\\
&+\frac{2r}{\varepsilon }\int_0^t\langle u^\varepsilon,\sin
u^\varepsilon\rangle_{L^2(D)}ds-\frac{2r}{\varepsilon
}\int_0^t\langle
\delta^\varepsilon,v_r^\varepsilon\rangle_{L^2(\partial D)}ds
+\frac{2r}{\varepsilon }\int_0^t\langle \theta_r^\varepsilon,
u^\varepsilon\rangle_{L^2(\partial D)}ds\\
&+\int_0^t\langle 2v_r^\varepsilon, \varepsilon^{\alpha-1}
dW_1(s)\rangle_{L^2(D)}+\int_0^t\langle 2\theta_r^\varepsilon,
\varepsilon^{\alpha-1} dW_2(s)\rangle_{L^2(\partial
D)}\\
&+\varepsilon^{2\alpha-2} TrQ_1\cdot t+\varepsilon^{2\alpha-2}
TrQ_2\cdot t.
\end{array}
\end{equation}
Meanwhile, we observe that
$$
\begin{array}{ll}
& \langle u^\varepsilon(t),\delta^\varepsilon(t)\rangle_{L^2(\partial D)}\\
= & \langle
u^\varepsilon(0),\delta^\varepsilon(0)\rangle_{L^2(\partial
D)}+\int_0^t \langle (u^\varepsilon)_s,\delta^\varepsilon
\rangle_{L^2(\partial D)}ds+\int_0^t \langle
u^\varepsilon,(\delta^\varepsilon)_s \rangle_{L^2(\partial D)}ds\\
= & \langle
u^\varepsilon(0),\delta^\varepsilon(0)\rangle_{L^2(\partial
D)}+\int_0^t \langle v^\varepsilon,\delta^\varepsilon
\rangle_{L^2(\partial D)}ds+\int_0^t \langle
u^\varepsilon,\theta^\varepsilon \rangle_{L^2(\partial D)}ds\\
= & \langle
u^\varepsilon(0),\delta^\varepsilon(0)\rangle_{L^2(\partial
D)}+\int_0^t \langle v_r^\varepsilon-r
u^\varepsilon,\delta^\varepsilon \rangle_{L^2(\partial
D)}ds+\int_0^t \langle
u^\varepsilon,\theta_r^\varepsilon-r\delta^\varepsilon \rangle_{L^2(\partial D)}ds\\
= & \langle
u^\varepsilon(0),\delta^\varepsilon(0)\rangle_{L^2(\partial
D)}+\int_0^t \langle v_r^\varepsilon,\delta^\varepsilon
\rangle_{L^2(\partial D)}ds+\int_0^t \langle
u^\varepsilon,\theta_r^\varepsilon\rangle_{L^2(\partial D)}ds\\
&-2r\int_0^t \langle
u^\varepsilon,\delta^\varepsilon\rangle_{L^2(\partial D)}ds,
\end{array}
$$
which implies that
\begin{equation}\label{Eq-3.1-8}
\begin{array}{ll}
&-\frac{2r}{\varepsilon }\int_0^t \langle
v_r^\varepsilon,\delta^\varepsilon \rangle_{L^2(\partial D)}ds
\\
=&-\frac{2r}{\varepsilon } \langle
u^\varepsilon(t),\delta^\varepsilon(t)\rangle_{L^2(\partial
D)}+\frac{2r}{\varepsilon } \langle
u^\varepsilon(0),\delta^\varepsilon(0)\rangle_{L^2(\partial
D)}+\frac{2r}{\varepsilon }\int_0^t \langle
u^\varepsilon,\theta_r^\varepsilon\rangle_{L^2(\partial D)}ds\\
&-\frac{4r^2}{\varepsilon }\int_0^t \langle
u^\varepsilon,\delta^\varepsilon\rangle_{L^2(\partial D)}ds.
\end{array}
\end{equation}
Hence, it   follows from (\ref{Eq-3.1-7}) and (\ref{Eq-3.1-8}) that
(\ref{Eq-3.1-01}) and (\ref{Eq-3.1-02}) hold. \hfill$\blacksquare$
\par
{\bf Proposition 3.2}\quad {\it Let $\alpha\in [1/2, 1)\cup
(1,+\infty)$. Assume that the initial datum $U_r^\varepsilon(0)$ is
a $ \mathcal{F}_0$-measurable random variable in $L^2(\Omega,
\mathcal{H} )$. Then for any time $t\in [0, \tau^*)$,
$\varepsilon\in (0,1/2)$  and a sufficiently small $r\in (0, 1/2)$,
there exists a positive constant $C$, independent of the parameter
$\varepsilon$, such that
\begin{equation}\label{Eq-3.2-0}
\begin{array}{ll}
&\frac{d}{dt}[\varepsilon
\mathbb{E}\|v_r^\varepsilon\|_{L^2(D)}^2+\mathbb{E}\|\bigtriangledown
u^\varepsilon\|_{L^2(D)}^2+\mathbb{E}\|u^\varepsilon\|_{L^2(D)}^2
+\varepsilon\mathbb{E}\|\theta_r^\varepsilon\|_{L^2(\partial
D)}^2+\mathbb{E}\|\delta^\varepsilon\|_{L^2(\partial D)}^2 ]\\
\leq &-C[\varepsilon
\mathbb{E}\|v_r^\varepsilon\|_{L^2(D)}^2+\mathbb{E}\|\bigtriangledown
u^\varepsilon\|_{L^2(D)}^2+\mathbb{E}\|u^\varepsilon\|_{L^2(D)}^2
+\varepsilon\mathbb{E}\|\theta_r^\varepsilon\|_{L^2(\partial
D)}^2+\mathbb{E}\|\delta^\varepsilon\|_{L^2(\partial
D)}^2]\\
&+C[TrQ_1 +TrQ_2+1].
\end{array}
\end{equation}
}
\par
{\bf Proof.}\quad On the one hand, it follows from the Cauchy
inequality and the trace inequality that there exists a positive
constant $C_{TI}>0$ (here and hereafter $C_{TI}$ denotes the
positive constant in the trace inequality) such that
$$
\begin{array}{ll}
0 &\leq r\mathbb{E}\|u^\varepsilon(t)\|_{L^2(\partial D)}^2+2r
\mathbb{E}\langle
u^\varepsilon(t),\delta^\varepsilon(t)\rangle_{L^2(\partial D)}+r
\mathbb{E}\|\delta^\varepsilon(t)\|_{L^2(\partial D)}^2\\
&\leq r C_{TI}^2\mathbb{E}\|u^\varepsilon(t)\|_{H ^1(D)}^2+2r
\mathbb{E}\langle
u^\varepsilon(t),\delta^\varepsilon(t)\rangle_{L^2(\partial D)}+r
\mathbb{E}\|\delta^\varepsilon(t)\|_{L^2(\partial D)}^2,
\end{array}
$$
which implies that
\begin{equation}\label{Eq-3.2-1}
\begin{array}{ll}
\mathbb{E} \mathcal{E}_r^\varepsilon(t)\geq &\varepsilon
\mathbb{E}\|v_r^\varepsilon(t)\|_{L^2(D)}^2+(1-r
C_{TI}^2)\mathbb{E}\|\bigtriangledown
u^\varepsilon(t)\|_{L^2(D)}^2\\
&+(1-r-r C_{TI}^2+\varepsilon
r^2)\mathbb{E}\|u^\varepsilon(t)\|_{L^2(D)}^2 +\varepsilon
\mathbb{E}\|\theta_r^\varepsilon(t)\|_{L^2(\partial D)}^2\\
&+(1-2r+\varepsilon r^2)\mathbb{E}\|\delta(t)\|_{L^2(\partial D)}^2.
\end{array}
\end{equation}
\par

On the other hand, by the H$\ddot{o}$lder inequality,
the Young inequality and the trace inequality, we obtain that
$$
\begin{array}{ll}
\mathbb{E}\langle
u^\varepsilon,\theta_r^\varepsilon\rangle_{L^2(\partial D)}&\leq
\mathbb{E}\|u^\varepsilon\|_{L^2(\partial D)}\cdot\mathbb{E}\|\theta_r^\varepsilon\|_{L^2(\partial D)}\\
&\leq r
\mathbb{E}\|u^\varepsilon\|_{L^2(\partial D)}^2+\frac{1}{4r}\mathbb{E}\|\theta_r^\varepsilon\|_{L^2(\partial D)}^2\\
&\leq r C_{TI}^2 \mathbb{E}\|u^\varepsilon\|_{H
^1(D)}^2+\frac{1}{4r}\mathbb{E}\|\theta_r^\varepsilon\|_{L^2(\partial
D)}^2,
\end{array}
$$
which implies that
\begin{equation}\label{Eq-3.2-2}
4r \mathbb{E}\langle
u^\varepsilon,\theta_r^\varepsilon\rangle_{L^2(\partial D)}\leq 4r^2
C_{TI}^2\mathbb{E}\|\bigtriangledown u^\varepsilon\|_{L^2(D)}^2+4r^2
C_{TI}^2\mathbb{E}\|u^\varepsilon\|_{L^2(D)}^2
+\mathbb{E}\|\theta_r^\varepsilon\|_{L^2(\partial D)}^2.
\end{equation}
\par
At the same time, it follows from the Cauchy inequality and the
trace inequality that
\begin{equation}\label{Eq-3.2-3}
\begin{array}{ll}
-4r^2 \mathbb{E}\langle
u^\varepsilon,\delta^\varepsilon\rangle_{L^2(\partial D)}&\leq 2r^2
\mathbb{E}\|u^\varepsilon\|_{L^2(\partial D)}^2+2r^2 \mathbb{E}\|\delta^\varepsilon\|_{L^2(\partial D)}^2\\
&\leq 2r^2C_{TI}^2 \mathbb{E}\|\bigtriangledown
u^\varepsilon\|_{L^2(D)}^2+2r^2C_{TI}^2
\mathbb{E}\|u^\varepsilon\|_{L^2(D)}^2+2r^2
\mathbb{E}\|\delta^\varepsilon\|_{L^2(\partial D)}^2.
\end{array}
\end{equation}
Also  the Cauchy inequality  leads to
\begin{equation}\label{Eq-3.2-4}
\begin{array}{ll}
2r\mathbb{E}\langle u^\varepsilon,\sin
u^\varepsilon\rangle_{L^2(D)}&\leq r
\mathbb{E}\|u^\varepsilon\|_{L^2(D)}^2+ r \mathbb{E}\|\sin
u^\varepsilon\|_{L^2(D)}^2\\
&\leq r \mathbb{E}\|u^\varepsilon\|_{L^2(D)}^2+ C.
\end{array}
\end{equation}
\par
Then it follows from Proposition 3.1 and
(\ref{Eq-3.2-2})-(\ref{Eq-3.2-4}) that
\begin{equation}\label{Eq-3.2-5}
\begin{array}{ll}
\mathbb{E} \mathcal{E}_r^\varepsilon(t)\leq& \mathbb{E}
\mathcal{E}_r^\varepsilon(0)-\int_0^t[2(1-\varepsilon
r)\mathbb{E}\|v^\varepsilon\|_{L^2(D)}^2 +2r(1-3r
C_{TI}^2)\mathbb{E}\|\bigtriangledown
u^\varepsilon\|_{L^2(D)}^2\\
&+r[1-2r-6r C_{TI}^2+2\varepsilon
r^2]\mathbb{E}\|u^\varepsilon\|_{L^2(D)}^2\\
&+(1-2\varepsilon r)\mathbb{E}\|\theta^\varepsilon\|_{L^2(\partial
D)}^2 +2r(1-2r+\varepsilon
r^2)\mathbb{E}\|\delta^\varepsilon\|_{L^2(\partial D)}^2]ds\\
&+\varepsilon^{2\alpha-1} TrQ_1\cdot
t+\varepsilon^{2\alpha-1}TrQ_2\cdot t+Ct.
\end{array}
\end{equation}
For $\varepsilon\in (0,1/2)$, choose $r$ in $(0, 1/2)$
sufficiently small such that
\begin{equation}\label{Eq-3.2-6}
\min\{1-3r C_{TI}^2, 1-r-r C_{TI}^2+\varepsilon r^2, 1-2r-6r
C_{TI}^2+2\varepsilon r^2, 1-2r+\varepsilon r^2
 \}>0.
\end{equation}
Then
\begin{equation}\label{Eq-3.2-7}
2(1-\varepsilon r)>\varepsilon,\quad \hbox{and}\quad 1-2\varepsilon
r>\varepsilon.
\end{equation}
Furthermore, noticing that $\alpha\in [1/2, 1)\bigcup
(1,+\infty)$, we have
\begin{equation}\label{Eq-3.2-8}
0<\varepsilon^{2\alpha-1}\leq 1.
\end{equation}
Therefore, from (\ref{Eq-3.2-1}), (\ref{Eq-3.2-5})-(\ref{Eq-3.2-8}),
there exists a positive constant $C$ independent of the parameter
$\varepsilon$ such that (\ref{Eq-3.2-0}) holds. \hfill$\blacksquare$
\par
{\bf Proposition 3.3}\quad {\it Let $\alpha\in [1/2, 1)\cup
(1,+\infty)$ and $\varepsilon\in (0,1/2)$. Assume that the initial
datum $U_r^\varepsilon(0)$ is a $ \mathcal{F}_0$-measurable random
variable in $L^2(\Omega, \mathcal{H} )$. Then there exists a
positive constant $C$ independent of the parameter $\varepsilon$
such that
\begin{equation}\label{Eq-3.3-0}
\varepsilon
\mathbb{E}\|v_r^\varepsilon\|_{L^2(D)}^2+\mathbb{E}\|\bigtriangledown
u^\varepsilon\|_{L^2(D)}^2+\mathbb{E}\|u^\varepsilon\|_{L^2(D)}^2
+\varepsilon\mathbb{E}\|\theta_r^\varepsilon\|_{L^2(\partial
D)}^2+\mathbb{E}\|\delta^\varepsilon\|_{L^2(\partial D)}^2 \leq C,
\forall t\in [0, \tau^*).
\end{equation}}
\par
Proposition 3.3 is easily deduced from the Gronwall inequality and
Proposition 3.2.
\par
{\bf Remark 3.1}\quad{\it From Frigeri \cite{Frigeri}, for
$r\in(0,1/2)$, $\mathbb{E}\|U_r^\varepsilon\|_{\mathcal{H} }^2\geq
\frac{1}{2}\mathbb{E}\|U^\varepsilon\|_{\mathcal{H} }^2$. Therefore,
if we obtain the almost sure boundedness of $U_r^\varepsilon$ in
$\mathcal{H}$, we naturally derive the almost sure boundedness of
$U^\varepsilon$ in $\mathcal{H}$. But from Proposition 3.3, since
the parameter $\varepsilon$ disturbs the system (\ref{Eq1}), we can
not use the pseudo energy method to directly derive the almost sure
boundedness, while this method is effective for wave equations
without   the singular parameter $\varepsilon$ (see
Chen and Zhang \cite{Chen2}). }
\par
{\bf Remark 3.2}\quad {\it Although Proposition 3.3 does not implies
the almost sure boundedness of $U^\varepsilon$ in $\mathcal{H}$, we
can obtain that under the condition of Proposition 3.3, there exists
a positive constant $C$ independent of the parameter $\varepsilon$
such that
\begin{equation}\label{Eq-3.3-1}
\mathbb{E}\|\bigtriangledown u^\varepsilon\|_{L^2(D)}^2\leq C,\quad
\mathbb{E}\|u^\varepsilon\|_{L^2(D)}^2\leq C,\quad
\mathbb{E}\|\delta^\varepsilon\|_{L^2(\partial D)}^2 \leq C, \quad
\forall t\in [0, \tau^*).
\end{equation}
}
\par
In the following, we will continue to derive the almost sure
boundedness of the solution for the Cauchy problem (\ref{Eq3}).
\par
{\bf Proposition 3.4}\quad {\it Let $\alpha\in [1/2, 1)\cup
(1,+\infty)$ and $\varepsilon\in (0,1/2)$. Assume that the initial
datum $U^\varepsilon(0)$ is a $ \mathcal{F}_0$-measurable random
variable in $L^2(\Omega, \mathcal{H} )$. Then there exists a
positive constant $C$ independent of the parameter $\varepsilon$
such that
\begin{equation}\label{Eq-3.4-0}
\mathbb{E}\|v^\varepsilon\|_{L^2(D)}^2 \leq C,\quad
\mathbb{E}\|\theta^\varepsilon\|_{L^2(\partial D)}^2 \leq C,\quad
\forall\; t\in [0, \tau^*).
\end{equation} }
\par
{\bf Proof.}\quad Set $ M(v^\varepsilon(t))=\int_D |v^\varepsilon
(t)|^2dx$. For the second equation of (\ref{Eq2}), from the
It$\hat{o}$ formula, we get
\begin{equation}\label{Eq-3.4-1}
\begin{array}{ll}
M(v^\varepsilon(t))=& M(v^\varepsilon(0))+\int_0^t\langle
M^{\prime}(v^\varepsilon),
\varepsilon^{\alpha-1}dW_1(s)\rangle_{L^2(D)}\\
&+ \int_0^t\langle
M^{\prime}(v^\varepsilon),(\frac{1}{\varepsilon}\bigtriangleup
u^\varepsilon-\frac{1}{\varepsilon}u^\varepsilon-\frac{1}{\varepsilon}v^\varepsilon+\frac{1}{\varepsilon}\sin
u^\varepsilon)\rangle_{L^2(D)}ds\\
&+\int_0^t\frac{1}{2}Tr[M^{\prime\prime}(v^\varepsilon)
(\varepsilon^{\alpha-1}Q_1)^{\frac{1}{2}}(\varepsilon^{\alpha-1}Q_1^{\frac{1}{2}})^*]ds,
\end{array}
\end{equation}
with  $M^\prime(v^\varepsilon)=2v^\varepsilon$ and
$M^{\prime\prime}(v^\varepsilon)=2\varphi$ for any $\varphi$ in
$L^2(D)$. Thus we deduce that
\begin{equation}\label{Eq-3.4-2}
\begin{array}{ll}
&\langle
M^{\prime}(v^\varepsilon),(\frac{1}{\varepsilon}\bigtriangleup
u^\varepsilon-\frac{1}{\varepsilon}u^\varepsilon-\frac{1}{\varepsilon}v^\varepsilon+\frac{1}{\varepsilon}\sin
u^\varepsilon)\rangle_{L^2(D)}\\
=& -\frac{d}{ds}[\frac{1}{\varepsilon}\|\bigtriangledown
u^\varepsilon \|_{L^2(D)}^2+\frac{1}{\varepsilon} \|u^\varepsilon
\|_{L^2(D)}^2+\frac{4}{\varepsilon}\|\cos
\frac{u^\varepsilon }{2}\|_{L^2(D)}^2]\\
&-\frac{2}{\varepsilon }\|v^\varepsilon \|_{L^2(D)}^2+2\langle
v^\varepsilon ,\frac{\partial u^\varepsilon }{\partial {\bf
n}}\rangle_{L^2(\partial D)}.
\end{array}
\end{equation}
It immediately follows from (\ref{Eq-3.4-1}) and (\ref{Eq-3.4-2})
that
\begin{equation}\label{Eq-3.4-3}
\begin{array}{ll}
&\|v^\varepsilon
(t)\|_{L^2(D)}^2+\frac{1}{\varepsilon}\|\bigtriangledown
u^\varepsilon (t)\|_{L^2(D)}^2+\frac{1}{\varepsilon}\|u^\varepsilon
(t)\|_{L^2(D)}^2+\frac{4}{\varepsilon}\|\cos
\frac{u^\varepsilon (t)}{2}\|_{L^2(D)}^2\\
=& \|v^\varepsilon
(0)\|_{L^2(D)}^2+\frac{1}{\varepsilon}\|\bigtriangledown
u^\varepsilon (0)\|_{L^2(D)}^2+\frac{1}{\varepsilon}\|u^\varepsilon
(0)\|_{L^2(D)}^2+\frac{4}{\varepsilon}\|\cos
\frac{u^\varepsilon (0)}{2}\|_{L^2(D)}^2\\
&-\frac{2}{\varepsilon}\int_0^t\|v^\varepsilon
\|_{L^2(D)}^2ds+\frac{2}{\varepsilon}\int_0^t\langle
v^\varepsilon ,\frac{\partial u^\varepsilon }{\partial {\bf n}}\rangle_{L^2(\partial D)}ds\\
&+\int_0^t\langle 2v^\varepsilon,
\varepsilon^{\alpha-1}dW_1(s)\rangle_{L^2(D)}+\varepsilon^{2\alpha-2}TrQ_1\cdot
t.
\end{array}
\end{equation}
\par
Secondly, noticing that the fourth equation of (\ref{Eq2}) and
putting $M(\theta^\varepsilon )=\int_{\partial
D}|\theta^\varepsilon| ^2dx$, using the It$\hat{o}$ formula, we have
\begin{equation}\label{Eq-3.4-4}
\begin{array}{lll}
M(\theta^\varepsilon(t))&=& M(\theta_0^\varepsilon)+\int_0^t\langle
M^{\prime}(\theta^\varepsilon), \varepsilon^{\alpha-1}
dW_2(s)\rangle_{L^2({\partial D})}\\
&&+ \int_0^t\langle
M^{\prime}(\theta^\varepsilon),(-\frac{1}{\varepsilon}\theta^\varepsilon-\frac{1}{\varepsilon}\delta^\varepsilon
-\frac{1}{\varepsilon}v^\varepsilon)\rangle_{L^2({\partial D})}ds\\
&&+\int_0^t\frac{1}{2}Tr[M^{\prime\prime}(\theta^\varepsilon)
(\varepsilon^{\alpha-1}Q_2^{\frac{1}{2}})(\varepsilon^{\alpha-1}Q_2^{\frac{1}{2}})^*]ds,
\end{array}
\end{equation}
with   $M^\prime(\theta^\varepsilon)=2\theta^\varepsilon$ and
$M^{\prime\prime}(\theta^\varepsilon)=2\phi$ for any $\phi$ in
$L^2({\partial D})$. After some further calculation, we   conclude that
\begin{equation}\label{Eq-3.4-5}
\begin{array}{ll}
&\langle
M^{\prime}(\theta^\varepsilon),(-\frac{1}{\varepsilon}\theta^\varepsilon-\frac{1}{\varepsilon}\delta^\varepsilon
-\frac{1}{\varepsilon}v^\varepsilon)\rangle_{L^2({\partial D})}\\
=&-\frac{1}{\varepsilon}\frac{d}{ds}\|\delta^\varepsilon\|_{L^2({\partial
D})}^2-\frac{2}{\varepsilon}\|\theta^\varepsilon\|_{L^2({\partial
D})}^2-\frac{2}{\varepsilon}\langle \frac{\partial
u^\varepsilon}{\partial {\bf n}},v^\varepsilon\rangle_{L^2({\partial
D})}.
\end{array}
\end{equation}
Thus, by (\ref{Eq-3.4-4}) and (\ref{Eq-3.4-5}),
\begin{equation}\label{Eq-3.4-6}
\begin{array}{ll}
&\|\theta^\varepsilon(t)\|_{L^2({\partial D})}^2+\frac{1}{\varepsilon}\|\delta^\varepsilon(t)\|_{L^2({\partial D})}^2\\
=&\|\theta^\varepsilon(0)\|_{L^2({\partial
D})}^2+\frac{1}{\varepsilon}\|\delta^\varepsilon(0)\|_{L^2({\partial
D})}^2
-\frac{2}{\varepsilon}\int_0^t\|\theta^\varepsilon\|_{L^2({\partial
D})}^2ds-\frac{2}{\varepsilon}\int_0^t\langle
\frac{\partial u^\varepsilon}{\partial {\bf n}},v^\varepsilon\rangle_{L^2({\partial D})}ds\\
&+\int_0^t\langle 2\theta^\varepsilon,
\varepsilon^{\alpha-1}dW_2(s)\rangle_{L^2({\partial
D})}+\varepsilon^{2\alpha-2}TrQ_2\cdot t.
\end{array}
\end{equation}
Then it   follows from (\ref{Eq-3.4-3}) and
(\ref{Eq-3.4-6}) that
$$
\begin{array}{ll}
&\mathbb{E}\|v^\varepsilon
(t)\|_{L^2(D)}^2+\mathbb{E}\|\theta^\varepsilon(t)\|_{L^2({\partial
D})}^2+\frac{1}{\varepsilon}\mathbb{E}\|\bigtriangledown
u^\varepsilon
(t)\|_{L^2(D)}^2+\frac{1}{\varepsilon}\mathbb{E}\|u^\varepsilon
(t)\|_{L^2(D)}^2\\
&+\frac{1}{\varepsilon}\mathbb{E}\|\delta^\varepsilon(t)\|_{L^2({\partial
D})}^2+\frac{4}{\varepsilon}\mathbb{E}\|\cos
\frac{u^\varepsilon (t)}{2}\|_{L^2(D)}^2\\
=& \mathbb{E}\|v^\varepsilon
(0)\|_{L^2(D)}^2+\mathbb{E}\|\theta^\varepsilon(0)\|_{L^2({\partial
D})}^2+\frac{1}{\varepsilon}\mathbb{E}\|\bigtriangledown
u^\varepsilon
(0)\|_{L^2(D)}^2+\frac{1}{\varepsilon}\mathbb{E}\|u^\varepsilon
(0)\|_{L^2(D)}^2\\
&+\frac{1}{\varepsilon}\mathbb{E}\|\delta^\varepsilon(0)\|_{L^2({\partial
D})}^2+\frac{4}{\varepsilon}\mathbb{E}\|\cos \frac{u^\varepsilon
(0)}{2}\|_{L^2(D)}^2\\
&-\frac{2}{\varepsilon}\int_0^t[\mathbb{E}\|v^\varepsilon
\|_{L^2(D)}^2+\mathbb{E}\|\theta^\varepsilon\|_{L^2({\partial
D})}^2]ds\\
&+\varepsilon^{2\alpha-2}TrQ_1\cdot
t+\varepsilon^{2\alpha-2}TrQ_2\cdot t,
\end{array}
$$
which implies that
$$
\begin{array}{ll}
&\mathbb{E}\|v^\varepsilon
(t)\|_{L^2(D)}^2+\mathbb{E}\|\theta^\varepsilon(t)\|_{L^2({\partial
D})}^2\\
\leq & \mathbb{E}\|v^\varepsilon
(0)\|_{L^2(D)}^2+\mathbb{E}\|\theta^\varepsilon(0)\|_{L^2({\partial
D})}^2+\frac{1}{\varepsilon}\mathbb{E}\|\bigtriangledown
u^\varepsilon
(0)\|_{L^2(D)}^2+\frac{1}{\varepsilon}\mathbb{E}\|u^\varepsilon
(0)\|_{L^2(D)}^2\\
&+\frac{1}{\varepsilon}\mathbb{E}\|\delta^\varepsilon(0)\|_{L^2({\partial
D})}^2+\frac{4}{\varepsilon}\mathbb{E}\|\cos \frac{u^\varepsilon
(0)}{2}\|_{L^2(D)}^2\\
&-\frac{2}{\varepsilon}\int_0^t[\mathbb{E}\|v^\varepsilon
\|_{L^2(D)}^2+\mathbb{E}\|\theta^\varepsilon\|_{L^2({\partial
D})}^2]ds\\
&+\varepsilon^{2\alpha-2}TrQ_1\cdot
t+\varepsilon^{2\alpha-2}TrQ_2\cdot t.
\end{array}
$$
Therefore, we have
\begin{equation}\label{Eq-3.4-7}
\begin{array}{ll}
&\frac{d}{dt}[\mathbb{E}\|v^\varepsilon
(t)\|_{L^2(D)}^2+\mathbb{E}\|\theta^\varepsilon(t)\|_{L^2({\partial
D})}^2]\\
\leq &-\frac{2}{\varepsilon}[\mathbb{E}\|v^\varepsilon
\|_{L^2(D)}^2+\mathbb{E}\|\theta^\varepsilon\|_{L^2({\partial
D})}^2]+\varepsilon^{2\alpha-2}TrQ_1+\varepsilon^{2\alpha-2}TrQ_2.
\end{array}
\end{equation}
By the Gronwall inequality, and noticing that $\alpha\in [1/2,
1)\cup (1,+\infty)$ and $\varepsilon\in (0,1/2)$, it follows from
(\ref{Eq-3.4-7}) that for arbitrary $t\in [0, \tau^*)$,
\begin{equation}\label{Eq-3.4-8}
\begin{array}{ll}
&\mathbb{E}\|v^\varepsilon
(t)\|_{L^2(D)}^2+\mathbb{E}\|\theta^\varepsilon(t)\|_{L^2({\partial
D})}^2\\
\leq & [\mathbb{E}\|v^\varepsilon(0)
\|_{L^2(D)}^2+\mathbb{E}\|\theta^\varepsilon(0)\|_{L^2({\partial
D})}^2]e^{-\frac{2t}{\varepsilon}}
+\frac{1}{2}(\varepsilon^{2\alpha-1}TrQ_1+\varepsilon^{2\alpha-1}TrQ_2)\cdot(1-e^{-\frac{2t}{\varepsilon}})\\
\leq & [\mathbb{E}\|v^\varepsilon(0)
\|_{L^2(D)}^2+\mathbb{E}\|\theta^\varepsilon(0)\|_{L^2({\partial
D})}^2]+\frac{1}{2}(\varepsilon^{2\alpha-1}TrQ_1+\varepsilon^{2\alpha-1}TrQ_2)\\
\leq & [\mathbb{E}\|v^\varepsilon(0)
\|_{L^2(D)}^2+\mathbb{E}\|\theta^\varepsilon(0)\|_{L^2({\partial
D})}^2]+\frac{1}{2}(TrQ_1+TrQ_2).
\end{array}
\end{equation}
This completes the proof of Proposition 3.4.
\hfill$\blacksquare$
\par
{\bf Proposition 3.5}\quad {\it Let $\alpha\in [1/2, 1)\cup
(1,+\infty)$ and $\varepsilon\in (0,1/2)$. Assume that the initial
datum $U^\varepsilon(0)$ is a $ \mathcal{F}_0$-measurable random
variable in $L^2(\Omega, \mathcal{H} )$. Then the solution
$U^\varepsilon(t)$ of the Cauchy problem (\ref{Eq3}) globally exists
in $\mathcal{H} $, i.e. $\tau^*=+\infty$ almost surely. }
\par
From Proposition 2.1, Remark 3.2 and Proposition 3.4, using the
Borel-Cantelli lemma, we easily obtain Proposition 3.5.
\par
{\bf Remark 3.3 (Almost sure boundedness)}\quad {\it From Remark
3.2, Proposition 3.4 and Proposition 3.5, we know that the global
solution $U^\varepsilon(t)$ of the Cauchy problem (\ref{Eq3}) is
bounded in $\mathcal{H} $ almost surely.}
\par
Introduce another space
$$ \Sigma : =\{U^\varepsilon\in
H ^2(D)\times H ^1(D)\times H^{1/2}(\partial D)\times
H^{1/2}(\partial D)|\; \frac{\partial u^\varepsilon}{\partial {\bf
n}}=\theta^\varepsilon\; \hbox{on} \;
\partial D \}.$$
\par
{\bf Proposition 3.7} \quad{\it Let $\alpha\in [1/2, 1)\bigcup
(1,+\infty)$ and $\varepsilon\in (0,1/2)$. Assume that the initial
datum $U^\varepsilon(0)$ is a $ \mathcal{F}_0$-measurable random
variable in $L^2(\Omega, \Sigma )$. Then the global solution
$U^\varepsilon(t)$ of the Cauchy problem (\ref{Eq3}) is also bounded
in $\Sigma $ almost surely. }
\par
Using a similar process for proving the almost sure boundedness of
the solution for the Cauchy problem (\ref{Eq3}) in $\mathcal{H} $,
we can prove Proposition 3.7. It is omitted here.
\par
We now establish the tightness of solutions for the
system (\ref{Eq1}). To begin with it, we recall some related
results. Let $\mathcal{X}\subset \mathcal{Y}\subset \mathcal{Z}$ be
three reflective Banach spaces and $\mathcal{X}\hookrightarrow
\mathcal{Y}$ being a compact and dense embedding. Define a new Banach
space
$$
\mathcal{G}=\{\varphi: \varphi\in L^2(0,T; \mathcal{X}), \; \frac{d
\varphi}{dt}\in L^2(0, T; \mathcal{Z})\},
$$
with   norm
$$
\|\varphi\|_{\mathcal{G}}^2=\int_0^T\|\varphi(s)\|_{\mathcal{X}}^2ds+\int_0^T\|\frac{d
\varphi(s)}{ds}\|_{\mathcal{Z}}^2ds.
$$
\par
{\bf Lemma 3.1}$^{\cite{Lions}}$\quad{\it If $K$ is bounded in
$\mathcal{G}$, then $K$ is precompact in $L^2(0, T; \mathcal{Y})$.}
\par
{\bf Proposition 3.8 (Tightness)}\quad{\it Let $\alpha\in [1/2,
1)\cup (1,+\infty)$ and $\varepsilon\in (0,1/2)$. Assume that the
initial datum $U^\varepsilon(0)$ is a $ \mathcal{F}_0$-measurable
random variable in $L^2(\Omega, \mathcal{H} )$. Then for a
given positive $T$,   the solution
$u^\varepsilon(t)$ and $\delta^\varepsilon(t)$ of the system
(\ref{Eq1}) is tight in $L^2(0,T; L^2(D))$ and $L^2(0,T;
L^2(\partial D))$, respectively. }
\par
{\bf Proof.}\quad First, for the solution $u^\varepsilon$ of the
system (\ref{Eq1}), let $\mathcal{X}=H^1(D)$ and
$\mathcal{Y}=\mathcal{Z}=L^2(D)$. Then from Remark 3.3 and the
Chebyshev inequality, for any $\rho>0$, there exists a bounded ball
of radius $\rho$ centered at zero, $K_{\rho}\subset \mathcal{G}$, such
that $\mathbb{P}\{u^\varepsilon\in K_\rho\}>1-\rho$. By Lemma 3.1,
$K_\rho$ is precompact in $L^2(0, T; \mathcal{Y})$. Then the
solution $u^\varepsilon$ is tight in $L^2(0, T; L^2(D))$.
\par
Similarly, for $\delta^\varepsilon$, let $\mathcal{X}=\mathcal{Y}=\mathcal{Z}=L^2(\partial
D)$. Using the same process of $u^\varepsilon$, we see that
  $\delta^\varepsilon$ is tight in $L^2(0, T; L^2(\partial
D))$. \hfill$\blacksquare$

\renewcommand{\theequation}{\thesection.\arabic{equation}}
\setcounter{equation}{0}

\section{Approximating equation}

\quad\quad In this section,    we   use a splitting method \cite{WLR} to derive the
approximating equation of the system (\ref{Eq1}) for $\varepsilon$
sufficiently small, in the sense of probability distribution. We
consider the solutions of the system (\ref{Eq1}) in the weak sense.
The main result is as follows.
\par
{\bf Theorem 4.1 (Approximating equation)}\quad {\it For the system
(\ref{Eq1}), let the initial datum $(u_0, u_1, \delta_0,\delta_1)^T$
be a $ \mathcal{F}_0$-measurable random variable in $L^2(\Omega,
\mathcal{H} )$. Let $T$ be a given positive number. Then we have the
following conclusions:
\par
(i)\quad If $\alpha\in [1/2, 1)$, then for sufficiently small
$\varepsilon$,
\begin{equation}\label{Eq-4.1-01}
\begin{array}{l}
\|u^\varepsilon-\overline{u}^\varepsilon\|_{L^2(0,
T; L^2(D))}=O(\varepsilon^\alpha),\\
\|\delta^\varepsilon-\overline{\delta}^\varepsilon\|_{L^2(0, T;
L^2(\partial D))}=O(\varepsilon^\alpha),
\end{array}
\end{equation}
where $\overline{u}^\varepsilon$ and $\overline{\delta}^\varepsilon$
are the solutions of the following stochastic parabolic equation
with a dynamical boundary condition
\begin{equation}\label{Eq-4.1-02}
\left\{
\begin{array}{ll}
\overline{u}^\varepsilon_t-\bigtriangleup
\overline{u}^\varepsilon+\overline{u}^\varepsilon-\sin
\overline{u}^\varepsilon=\varepsilon^\alpha
\dot{\overline{W}}_1,\quad &in\quad D,\\
\overline{\delta}^\varepsilon_t+\overline{\delta}^\varepsilon =
-\overline{u}^\varepsilon_t+\varepsilon^\alpha
\dot{\overline{W}}_2,\quad &on\quad \partial D,\\
\overline{\delta}^\varepsilon_t =\frac{\partial
\overline{u}^\varepsilon}{\partial {\bf
n}},\quad &on\quad \partial D,\\
\overline{u}^\varepsilon (0)=u_0,
\overline{\delta}^\varepsilon(0)=\delta_0.
\end{array}
\right.
\end{equation}

(ii)\quad If $\alpha\in (1,+\infty)$, then for sufficiently small
$\varepsilon$,
\begin{equation}\label{Eq-4.1-03}
\begin{array}{l}
\|u^\varepsilon-\overline{u}^\varepsilon\|_{L^2(0,
T; L^2(D))}= O(\varepsilon),\\
\|\delta^\varepsilon-\overline{\delta}^\varepsilon\|_{L^2(0, T;
L^2(\partial D))}= O(\varepsilon),
\end{array}
\end{equation}
where $\overline{u}^\varepsilon$ and $\overline{\delta}^\varepsilon$
are the solutions of the following deterministic wave equation with
a dynamical boundary condition
\begin{equation}\label{Eq-4.1-04}
\left\{
\begin{array}{ll}
\varepsilon
\overline{u}^\varepsilon_{tt}+\overline{u}^\varepsilon_t-\bigtriangleup
\overline{u}^\varepsilon+\overline{u}^\varepsilon-\sin
\overline{u}^\varepsilon=0, &\quad in\; D,\\
\varepsilon\overline{\delta}^\varepsilon_{tt}+\overline{\delta}^\varepsilon_t+\overline{\delta}^\varepsilon
=-\overline{u}^\varepsilon_t, &\quad on\;
\partial D ,\\
\overline{\delta}^\varepsilon_t=\frac{\partial
\overline{u}^\varepsilon}{\partial {\bf n}}, &\quad on\;
\partial D,\\
\overline{u}^\varepsilon(0)=u_0, \overline{u}^\varepsilon_t(0)=u_1,
\overline{\delta}^\varepsilon(0)=\delta_0,
\overline{\delta}^\varepsilon_t(0)=\delta_1.
\end{array}
\right.
\end{equation}}
\par
{\bf Remark 4.1}\quad{\it By the method of Chueshov and
Schmalfuss \cite{CS04, CS07}, we can show that Equation (\ref{Eq-4.1-02}) is  well-posed. In addition, Equation (\ref{Eq-4.1-04})
  is also well-posed (see \cite{Frigeri}). }
\par
In the following, we will prove Theorem 4.1.  We
state some preliminary results.
\par
For the system (\ref{Eq2}), we give the decomposition as follows.
Firstly, \begin{equation}\label{Eq-4-1}
v^\varepsilon=\overline{v}_1^\varepsilon+\overline{v}_2^\varepsilon+\overline{v}_3^\varepsilon,
\end{equation}
where
\begin{equation}\label{Eq-4-1.1}
\left\{
\begin{array}{l}
\frac{d
\overline{v}_1^\varepsilon}{dt}=-\frac{1}{\varepsilon}\overline{v}_1^\varepsilon, \quad\hbox{in}\quad D,\\
\overline{v}_1^\varepsilon(0)=v_0,
\end{array}
\right.
\end{equation}
\begin{equation}\label{Eq-4-1.2}
\left\{
\begin{array}{l}
\frac{d
\overline{v}_2^\varepsilon}{dt}=-\frac{1}{\varepsilon}\overline{v}_2^\varepsilon
+\frac{1}{\varepsilon}\bigtriangleup u^\varepsilon-\frac{1}{\varepsilon}u^\varepsilon
+\frac{1}{\varepsilon}\sin u^\varepsilon, \quad\hbox{in}\quad D,\\
\overline{v}_2^\varepsilon(0)=0,
\end{array}
\right.
\end{equation}
and
\begin{equation}\label{Eq-4-1.3}
\left\{
\begin{array}{l}
\frac{d
\overline{v}_3^\varepsilon}{dt}=-\frac{1}{\varepsilon}\overline{v}_3^\varepsilon
+\varepsilon^{\alpha-1}\dot{W}_1, \quad\hbox{in}\quad D,\\
\overline{v}_3^\varepsilon(0)=0.
\end{array}
\right.
\end{equation}
\par
Secondly,
\begin{equation}\label{Eq-4-2}
\theta^\varepsilon=\overline{\theta}_1^\varepsilon+\overline{\theta}_2^\varepsilon+\overline{\theta}_3^\varepsilon,
\end{equation}
where
\begin{equation}\label{Eq-4-2.1}
\left\{
\begin{array}{l}
\frac{d
\overline{\theta}_1^\varepsilon}{dt}=-\frac{1}{\varepsilon}\overline{\theta}_1^\varepsilon,
\quad\hbox{on}\quad \partial D,\\
\overline{\theta}_1^\varepsilon(0)=\theta_0,
\end{array}
\right.
\end{equation}
\begin{equation}\label{Eq-4-2.2}
\left\{
\begin{array}{l}
\frac{d
\overline{\theta}_2^\varepsilon}{dt}=-\frac{1}{\varepsilon}\overline{\theta}_2^\varepsilon
-\frac{1}{\varepsilon}
\delta^\varepsilon+\frac{1}{\varepsilon}v^\varepsilon, \quad\hbox{on}\quad \partial D,\\
\overline{\theta}_2^\varepsilon(0)=0,
\end{array}
\right.
\end{equation}
and
\begin{equation}\label{Eq-4-2.3}
\left\{
\begin{array}{l}
\frac{d
\overline{\theta}_3^\varepsilon}{dt}=-\frac{1}{\varepsilon}\overline{\theta}_3^\varepsilon
+\varepsilon^{\alpha-1}\dot{W}_2, \quad\hbox{on}\quad \partial D,\\
\overline{\theta}_3^\varepsilon(0)=0.
\end{array}
\right.
\end{equation}
\par
{\bf Proposition 4.1}\quad{\it Let $\alpha\in [1/2, 1)\cup
(1,+\infty)$ and $\varepsilon\in (0,1/2)$. Assume that the initial
data $v_0$ and $\theta_0$ are $ \mathcal{F}_0$-measurable random
variables in $L^2(\Omega, L^2(D))$ and $L^2(\Omega, L^2(\partial
D))$, respectively. Then we have that
\par
(i)\quad For Equation (\ref{Eq-4-1.1}) and Equation (\ref{Eq-4-2.1})
\begin{equation}\label{Eq-p4.1-01}
\overline{v}_1^{\varepsilon}(t)=v_0e^{-\frac{t}{\varepsilon}}, \quad
\overline{\theta}_1^{\varepsilon}(t)=\theta_0e^{-\frac{t}{\varepsilon}},
\quad \forall\;t\geq 0.
\end{equation}
\par
(ii)\quad For Equation (\ref{Eq-4-1.2}) and Equation
(\ref{Eq-4-2.2}), there is a positive constant $C$ independent of
the parameter $\varepsilon$ such that
\begin{equation}\label{Eq-p4.1-02}
\mathbb{E}\|\overline{v}^\varepsilon_2(t)\|_{H^{-1}(D)}\leq C,\quad
\mathbb{E}\|\overline{\theta}^\varepsilon_2(t)\|_{H^{-1/2}(\partial
D)}\leq C, \quad\forall \; t\geq 0.
\end{equation}
\par
(iii)\quad For Equation (\ref{Eq-4-1.3}) and Equation
(\ref{Eq-4-2.3}),
\begin{equation}\label{Eq-p4.1-03}
\begin{array}{l}
\mathbb{E}\|\overline{v}_3^{\varepsilon}(t)\|_{L^2(D)}^2
=-\frac{2}{\varepsilon}\int_0^t\mathbb{E}\|\overline{v}_3^{\varepsilon}(s)\|_{L^2(D)}^2ds
+\varepsilon^{2\alpha-2}TrQ_1\cdot t,\quad \forall\;t\geq 0,\\
\mathbb{E}\|\overline{\theta}_3^{\varepsilon}(t)\|_{L^2(D)}^2
=-\frac{2}{\varepsilon}\int_0^t\mathbb{E}\|\overline{\theta}_3^{\varepsilon}(s)\|_{L^2(\partial
D)}^2ds +\varepsilon^{2\alpha-2}TrQ_2\cdot t,\quad \forall\;t\geq 0.
\end{array}
\end{equation}
}
\par
{\bf Proof.}\quad For Equation (\ref{Eq-4-1.1}) and Equation
(\ref{Eq-4-2.1}), we can directly solve them to obtain
(\ref{Eq-p4.1-01}). And for Equation (\ref{Eq-4-1.3}) and Equation
(\ref{Eq-4-2.3}), applying the It$\hat{o}$ formula, we immediately
obtain (\ref{Eq-p4.1-03}).
\par
Now we   prove (\ref{Eq-p4.1-02}).
\par
Noticing that
$$
\begin{array}{l}
\mathbb{E}\|\overline{v}^\varepsilon_2(t)\|_{H^{-1}(D)}=\sup\limits_{\phi\in
H^1(D)}\frac{\mathbb{E}|\langle\overline{v}^\varepsilon_2,
\phi\rangle_{L^2(D)}|}{\|\phi\|_{H^1(D)}},\\
\mathbb{E}\|\overline{\theta}^\varepsilon_2(t)\|_{H^{-1/2}(\partial
D)}=\sup\limits_{\psi\in H^{1/2}(\partial
D)}\frac{\mathbb{E}|\langle\overline{\theta}^\varepsilon_2,
\psi\rangle_{L^2(\partial D)}|}{\|\psi\|_{H^{1/2}(\partial D)}},
\end{array}
$$
we only need to prove
\begin{equation}\label{Eq-p4.1-1}
\mathbb{E}|\langle\overline{v}^\varepsilon_2,
\phi\rangle_{L^2(D)}|\leq C\|\phi\|_{H^1(D)},\quad \forall\;\phi\in
H^1(D).
\end{equation}
and
\begin{equation}\label{Eq-p4.1-2}
\mathbb{E}|\langle\overline{\theta}^\varepsilon_2,
\psi\rangle_{L^2(\partial D)}|\leq C\|\psi\|_{H^{1/2}(D)},\quad
\forall\;\psi\in H^{1/2}(D).
\end{equation}
\par
Firstly, for arbitrary $\phi\in H^1(D)$, it follows from
(\ref{Eq-4-1.2}) that
$$
\frac{d}{dt}\langle\overline{v}_2^\varepsilon,
\phi\rangle_{L^2(D)}=-\frac{1}{\varepsilon}\langle
\overline{v}_2^\varepsilon, \phi\rangle_{L^2(D)}
+\frac{1}{\varepsilon}\langle \bigtriangleup u^\varepsilon,
\phi\rangle_{L^2(D)}-\frac{1}{\varepsilon}\langle u^\varepsilon,
\phi\rangle_{L^2(D)} +\frac{1}{\varepsilon}\langle\sin
u^\varepsilon, \phi\rangle_{L^2(D)},
$$
which implies, from Remark 3.3, that
$$
\begin{array}{lll}
\mathbb{E}\langle\overline{v}_2^\varepsilon,
\phi\rangle_{L^2(D)}&=&\frac{1}{\varepsilon}
e^{-\frac{t}{\varepsilon}}\int_0^te^{\frac{s}{\varepsilon}}
[-\mathbb{E}\langle \bigtriangledown u^\varepsilon,
\bigtriangledown\phi\rangle_{L^2(D)}+\mathbb{E}\langle
\frac{\partial u^\varepsilon}{\partial {\bf n}},
\phi\rangle_{L^2(\partial D)}- \mathbb{E}\langle
u^\varepsilon, \phi\rangle_{L^2(D)}\\
&& +\mathbb{E}\langle\sin
u^\varepsilon, \phi\rangle_{L^2(D)}]ds\\
& \leq& \frac{1}{\varepsilon} e^{-\frac{t}{\varepsilon}}\int_0^t
e^{\frac{s}{\varepsilon}} [\mathbb{E}\| \bigtriangledown
u^\varepsilon\|_{L^2(D)}\cdot\|\phi\|_{H^1(D)}+\mathbb{E}\|
\theta^\varepsilon\|_{L^2(\partial D)}\cdot\|\phi\|_{H^1(D)}\\
&&+\mathbb{E}\| u^\varepsilon\|_{L^2(D)}\cdot
\|\phi\|_{H^1(D)}+\mathbb{E}\|\sin
u^\varepsilon\|_{L^2(D)}\cdot\|\phi\|_{H^1(D)}]ds\\
&\leq & \frac{1}{\varepsilon} e^{-\frac{t}{\varepsilon}}\int_0^t
e^{\frac{s}{\varepsilon}}ds\cdot C \|\phi\|_{H^1(D)}\\
&= & [1-e^{-\frac{t}{\varepsilon}}]\cdot C \|\phi\|_{H^1(D)}\\
&\leq & C \|\phi\|_{H^1(D)}.
\end{array}
$$
which arrives at (\ref{Eq-p4.1-1}).
\par
Secondly, for arbitrary $\psi\in H^{1/2}(\partial D)$, it follows
from (\ref{Eq-4-2.2}) that
$$
\frac{d}{dt}\langle\overline{\theta}_2^\varepsilon,
\psi\rangle_{L^2(\partial D)}=-\frac{1}{\varepsilon}\langle
\overline{\theta}_2^\varepsilon, \psi\rangle_{L^2(\partial D)}
-\frac{1}{\varepsilon}\langle \delta^\varepsilon,
\psi\rangle_{L^2(\partial D)}+\frac{1}{\varepsilon}\langle
v^\varepsilon, \psi\rangle_{L^2(\partial D)},
$$
which implies, from the trace inequality, Remark 3.3 and Proposition
3.7, that
$$
\begin{array}{lll}
\mathbb{E}\langle\overline{\theta}_2^\varepsilon,
\psi\rangle_{L^2(\partial D)}&=&\frac{1}{\varepsilon}
e^{-\frac{t}{\varepsilon}}\int_0^te^{\frac{s}{\varepsilon}}
[-\mathbb{E}\langle \delta^\varepsilon, \psi\rangle_{L^2(\partial
D)}+\mathbb{E}\langle
v^\varepsilon, \psi\rangle_{L^2(\partial D)}]ds\\
& \leq& \frac{1}{\varepsilon} e^{-\frac{t}{\varepsilon}}\int_0^t
e^{\frac{s}{\varepsilon}}
[\mathbb{E}\|\delta^\varepsilon\|_{L^2(\partial D)}^2\cdot
\|\psi\|_{H^{1/2}(\partial D)}+C_{TI}\mathbb{E}\|
v^\varepsilon\|_{H^1(D)}\cdot\|\psi\|_{H^{1/2}(\partial D)}]ds\\
&\leq & \frac{1}{\varepsilon} e^{-\frac{t}{\varepsilon}}\int_0^t
e^{\frac{s}{\varepsilon}}ds\cdot C \|\psi\|_{H^{1/2}(\partial D)}\\
&= &[1-e^{-\frac{t}{\varepsilon}}]\cdot C \|\psi\|_{H^{1/2}(\partial
D)}\\
&\leq & C \|\psi\|_{H^{1/2}(\partial D)}.
\end{array}
$$
which leads to (\ref{Eq-p4.1-2}). \hfill$\blacksquare$
\par
{\bf Lemma 4.1 (Prohorov Theorem)}$^{\cite{DZ1}}$\quad {\it Assume that
$\mathcal{M}$ is a separable Banach space. The set of probability
measures $\{\mathcal{L}(X_n)\}_n$ on $(\mathcal{M},
\mathcal{B}(\mathcal{M}))$ is relatively compact if and only if
$\{X_n\}$ is tight.}
\par
{\bf Lemma 4.2 (Skorohod Theorem)}$^{\cite{DZ1}}$\quad {\it For an
arbitrary sequence of probability measures $\{\mu_n\}$ on
$(\mathcal{M}, \mathcal{B}(\mathcal{M}))$ weakly converges to
probability measures $\mu$, there exists a probability space
$(\Omega, \mathcal{F}, \mathbb{P})$ and random variables, $X$,
$X_1$, $X_2$, $\cdots$, $X_n$, $\cdots$ such that $X_n$ distributes
as $\mu_n$ and $X$ distributes as $\mu$, and $\lim\limits_{n\to
\infty}X_n=X$, $\mathbb{P}$-a.s.}
\par
\vspace{1cm}
\par
{\bf Proof of Theorem 4.1}
\par

From Proposition 3.8, for   $t\in [0,
T]$, the solution $u^\varepsilon(t)$ and $\delta^\varepsilon(t)$ of
the system (\ref{Eq1}), are tight in $L^2(0,T; L^2(D))$ and
$L^2(0,T; L^2(\partial D))$, respectively. Therefore, for arbitrary
$\rho>0$, there exist two bounded balls of radius $\rho$ centered at
zero, $K_\rho\subset H^1(D)$ and $B_\rho\subset L^2(\partial D)$,
which are compact in $L^2(D)$ and $L^2(\partial D)$, such that
$$
\mathbb{P}\{u^\varepsilon\in K_\rho\}>1-\rho, \quad \hbox{and}\quad
\mathbb{P}\{\delta^\varepsilon\in B_\rho\}>1-\rho.
$$
According to Lemma 4.1 and Lemma 4.2, we know that for every sequence
$\{\varepsilon_j\}_{j=1}^{j=\infty}$ with $\varepsilon_j\to 0$ as
$j\to \infty$, there exists a subsequence
$\{\varepsilon_{j(k)}\}_{k=1}^{k=\infty}$, random variables
$u^{*\varepsilon_{j(k)}}\subset L^2(0,T;L^2(D))$ and
$\delta^{*\varepsilon_{j(k)}}\subset L^2(0,T;L^2(\partial D))$, and
$u^*\in L^2(0,T;L^2(D))$ and $\delta^*\in L^2(0,T;L^2(\partial D))$
defined on a new probability space $(\Omega^*, \mathcal{F}^*,
\mathbb{P}^*)$, such that for almost all $\omega\in \Omega^*$,
$$
\begin{array}{l}
\mathcal{L}(u^{*\varepsilon_{j(k)}})=\mathcal{L}(u^{\varepsilon_{j(k)}}),\\
\mathcal{L}(\delta^{*\varepsilon_{j(k)}})=\mathcal{L}(\delta^{\varepsilon_{j(k)}}),
\end{array}
$$
and
$$
\begin{array}{l}
u^{*\varepsilon_{j(k)}} \longrightarrow  u^*, \quad \hbox{in}\quad
L^2(0,T; L^2(D)) \quad \hbox{as}\quad k\to \infty,\\
\delta^{*\varepsilon_{j(k)}} \longrightarrow  \delta^*, \quad
\hbox{in}\quad L^2(0,T; L^2(\partial D)) \quad \hbox{as}\quad k\to
\infty.
\end{array}
$$
In the meantime, $u^{*\varepsilon_{j(k)}}$ and
$\delta^{*\varepsilon_{j(k)}}$ solve the system (\ref{Eq1}) with
$W_1$ and $W_2$ being replaced by the Wiener processes $W_1^*$ and
$W_2^*$, defined on the probability space $(\Omega^*, \mathcal{F}^*,
\mathbb{P}^*)$ but with the same distributions as $W_1$ and $W_2$,
respectively. In the following, we will derive the approximating
equation for $u^*$ and $\delta^*$ and present the error estimates between the
approximating equation and the original system (\ref{Eq1}) as in Theorem 4.1.
\par
Now, for the above $\rho$, it follows from (\ref{Eq-p4.1-02}) and
the Chebyshev inequality that there exists a positive constant $C_\rho$
independent of the parameter $\varepsilon$ such that
$$
\mathbb{P}\{\|\overline{v}_2^\varepsilon\|_{H^{-1}(D)}\leq
C_\rho\}>1-\rho, \quad \hbox{and}\quad
\mathbb{P}\{\|\overline{\theta}_2^\varepsilon\|_{H^{-1/2}(\partial
D)}\leq C_\rho\}>1-\rho.
$$
\par
Define
$$
\Omega_\rho=\{\omega\in\Omega: u^\varepsilon(\omega)\in K_\rho, \;
\delta^\varepsilon(\omega)\in B_\rho,\;
\|\overline{v}_2^\varepsilon(\omega)\|_{H^{-1}(D)}\leq C_\rho,\;
\|\overline{\theta}_2^\varepsilon(\omega)\|_{H^{-1/2}(\partial
D)}\leq C_\rho \},
$$
$$
\mathcal{F}_\rho=\{F\bigcap\Omega_\rho: F\in \mathcal{F}\},
$$
and
$$
\mathbb{P}_\rho(F)=\frac{\mathbb{P}(F\bigcap\Omega_\rho)}{\mathbb{P}(\Omega_\rho)},
\quad \hbox{for}\quad F\in \mathcal{F}_\rho.
$$
Then $(\Omega_\rho, \mathcal{F}_\rho, \mathbb{P}_\rho)$ is  a new
probability space, whose expectation operator is denoted by
$\mathbb{E}_\rho$. We will work in the probability
space $(\Omega_\rho, \mathcal{F}_\rho, \mathbb{P}_\rho)$ in stead of
$(\Omega, \mathcal{F}, \mathbb{P})$. For simplicity, we will omit
the subscript $\rho$ unless we specifically stated otherwise.
\par
The system (\ref{Eq2}), combining with (\ref{Eq-4-1}) and
(\ref{Eq-4-2}),   can be rewritten as follows
$$
\left\{
\begin{array}{ll}
u_t^\varepsilon=v^\varepsilon=\overline{v}_1^\varepsilon
+\overline{v}_2^\varepsilon+\overline{v}_3^\varepsilon,\quad u^\varepsilon(0)=u_0,&\quad \hbox{in}\quad D,\\
\delta_t^\varepsilon=\theta^\varepsilon=\overline{\theta}_1^\varepsilon
+\overline{\theta}_2^\varepsilon+\overline{\theta}_3^\varepsilon,\quad
\delta^\varepsilon(0)=\delta_0,&\quad \hbox{on}\quad \partial D,\\
\delta_t^\varepsilon=\frac{\partial u^\varepsilon}{\partial {\bf
n}},&\quad \hbox{on}\quad \partial D,
\end{array}
\right.
$$
whose weak sense formulation is
\begin{equation}\label{Eq-4.1-1}
\begin{array}{ll}
&\langle u^\varepsilon(t), \varphi(t) \rangle_{L^2(D)} +\langle
\delta^\varepsilon(t), \varphi(t) \rangle_{L^2(\partial
D)}\\
= &\langle u^\varepsilon(0), \varphi(0) \rangle_{L^2(D)}+\int_0^t
\langle
u^\varepsilon(s), \varphi_s(s) \rangle_{L^2(D)}ds \\
&+\langle \delta^\varepsilon(0), \varphi(0) \rangle_{L^2(\partial
D)}+\int_0^t \langle \delta^\varepsilon(s), \varphi_s(s)
\rangle_{L^2(\partial D)}ds\\
&+ \int_0^t\langle \overline{v}_1^\varepsilon(s), \varphi(s)
\rangle_{L^2(D)}ds + \int_0^t\langle \overline{v}_2^\varepsilon(s),
\varphi(s) \rangle_{L^2(D)}ds + \int_0^t\langle
\overline{v}_3^\varepsilon(s), \varphi(s) \rangle_{L^2(D)}ds\\
& + \int_0^t\langle \overline{\theta}_1^\varepsilon(s), \varphi(s)
\rangle_{L^2(\partial D)}ds + \int_0^t\langle
\overline{\theta}_2^\varepsilon(s), \varphi(s) \rangle_{L^2(\partial
D)}ds + \int_0^t\langle \overline{\theta}_3^\varepsilon(s),
\varphi(s) \rangle_{L^2(\partial D)}ds,
\end{array}
\end{equation}
for every $\varphi\in C_0^\infty([0,+\infty)\times D)$.
\par

We   consider the case of $\alpha\in [1/2,
1)$.
\par
From (\ref{Eq-p4.1-01}), it follows that for every $\varphi\in
C_0^\infty([0,+\infty)\times D)$,
\begin{equation}\label{Eq-4.1-2}
\begin{array}{ll}
\int_0^t\langle \overline{v}_1^{\varepsilon},
\varphi\rangle_{L^2(D)}ds &= \int_0^t\langle
v_0e^{-\frac{s}{\varepsilon}}, \varphi(s)\rangle_{L^2(D)}ds\\
&= \varepsilon \int_0^{\frac{t}{\varepsilon}} \langle v_0,
\varphi(\varepsilon \tau)\rangle_{L^2(D)} e^{-\tau}d\tau\\
&=O(\varepsilon),
\end{array}
\end{equation}
for sufficiently small $\varepsilon$.
\par
In addition, it follows from (\ref{Eq-4-1.2}) that
\begin{equation}\label{Eq-4.1-3}
\begin{array}{ll}
\int_0^t\langle \frac{d\overline{v}_2^{\varepsilon}}{ds},
\varphi\rangle_{L^2(D)}ds =&-\frac{1}{\varepsilon} \int_0^t\langle
\overline{v}_2^{\varepsilon},
\varphi\rangle_{L^2(D)}ds+\frac{1}{\varepsilon} \int_0^t\langle
\bigtriangleup u^{\varepsilon},
\varphi\rangle_{L^2(D)}ds\\
&-\frac{1}{\varepsilon} \int_0^t\langle u^{\varepsilon},
\varphi\rangle_{L^2(D)}ds+\frac{1}{\varepsilon} \int_0^t\langle \sin
u^{\varepsilon}, \varphi\rangle_{L^2(D)}ds.
\end{array}
\end{equation}
Meanwhile, noticing that $\overline{v}_2^{\varepsilon}(0)=0$
and that
$$
\int_0^t\langle \frac{d\overline{v}_2^{\varepsilon}}{ds},
\varphi\rangle_{L^2(D)}ds= \langle \overline{v}_2^{\varepsilon}(t),
\varphi(t)\rangle_{L^2(D)}-\langle \overline{v}_2^{\varepsilon}(0),
\varphi(0)\rangle_{L^2(D)}-\int_0^t\langle
\overline{v}_2^{\varepsilon}(s), \varphi_s(s)\rangle_{L^2(D)}ds,
$$
we infer from (\ref{Eq-4.1-3}) that
\begin{equation}\label{Eq-4.1-4}
\begin{array}{ll}
&\int_0^t \langle \overline{v}_2^{\varepsilon},
\varphi\rangle_{L^2(D)}ds\\
=&-\varepsilon \langle
\overline{v}_2^{\varepsilon}(t), \varphi(t)\rangle_{L^2(D)}+
\varepsilon \int_0^t\langle \overline{v}_2^{\varepsilon}(s),
\varphi_s(s)\rangle_{L^2(D)}ds\\
&+\int_0^t\langle \bigtriangleup u^{\varepsilon},
\varphi\rangle_{L^2(D)}ds- \int_0^t\langle u^{\varepsilon},
\varphi\rangle_{L^2(D)}ds+\int_0^t\langle \sin u^{\varepsilon},
\varphi\rangle_{L^2(D)}ds,
\end{array}
\end{equation}
which implies from (\ref{Eq-p4.1-02}) that
\begin{equation}\label{Eq-4.1-5}
\begin{array}{ll}
&\int_0^t \langle \overline{v}_2^{\varepsilon},
\varphi\rangle_{L^2(D)}ds\\
=& O(\varepsilon)+O(\varepsilon)+\int_0^t\langle \bigtriangleup
u^{\varepsilon}, \varphi\rangle_{L^2(D)}ds- \int_0^t\langle
u^{\varepsilon}, \varphi\rangle_{L^2(D)}ds+\int_0^t\langle \sin
u^{\varepsilon}, \varphi\rangle_{L^2(D)}ds,
\end{array}
\end{equation}
for sufficiently small $\varepsilon$.
\par
Also, it follows from (\ref{Eq-4-1.3}) that
\begin{equation}\label{Eq-4.1-6}
\int_0^t\langle \frac{d\overline{v}_3^{\varepsilon}}{ds},
\varphi\rangle_{L^2(D)}ds =-\frac{1}{\varepsilon} \int_0^t\langle
\overline{v}_3^{\varepsilon},
\varphi\rangle_{L^2(D)}ds+\varepsilon^{\alpha-1}\int_0^t \langle
dW_1(s), \varphi\rangle_{L^2(D)}.
\end{equation}
Noticing that $\overline{v}_3^{\varepsilon}(0)=0$ and that
$$
\int_0^t\langle \frac{d\overline{v}_3^{\varepsilon}}{ds},
\varphi\rangle_{L^2(D)}ds= \langle \overline{v}_3^{\varepsilon}(t),
\varphi(t)\rangle_{L^2(D)}-\langle \overline{v}_3^{\varepsilon}(0),
\varphi(0)\rangle_{L^2(D)}-\int_0^t\langle
\overline{v}_3^{\varepsilon}(s), \varphi_s(s)\rangle_{L^2(D)}ds,
$$
we deduce from (\ref{Eq-4.1-6}) that
\begin{equation}\label{Eq-4.1-7}
\begin{array}{ll}
\int_0^t\langle \overline{v}_3^{\varepsilon},
\varphi\rangle_{L^2(D)}ds = &-\varepsilon \langle
\overline{v}_3^{\varepsilon}(t),
\varphi(t)\rangle_{L^2(D)}ds+\varepsilon\int_0^t\langle
\overline{v}_3^{\varepsilon}(s),
\varphi_s(s)\rangle_{L^2(D)}ds\\
&+\varepsilon^\alpha\int_0^t \langle
dW_1(s), \varphi\rangle_{L^2(D)}.
\end{array}
\end{equation}
Combining with (\ref{Eq-p4.1-03}) and $\alpha\in [1/2, 1)\cup
(1,+\infty)$, and using the Gronwall inequality, we further obtain that
\begin{equation}\label{Eq-4.1-8}
\mathbb{E}\|\overline{v}_3^{\varepsilon}(t)\|_{L^2(D)}\leq
TrQ_1,\quad \forall\; t\geq 0,
\end{equation}
which immediately implies from (\ref{Eq-4.1-7}) and $\alpha\in [1/2,
1)$ that,
\begin{equation}\label{Eq-4.1-9}
\begin{array}{l}
\int_0^t\langle \overline{v}_3^{\varepsilon},
\varphi\rangle_{L^2(D)}ds =
O(\varepsilon)+O(\varepsilon)+\varepsilon^\alpha\int_0^t \langle
dW_1(s), \varphi\rangle_{L^2(D)},
\end{array}
\end{equation}
for sufficiently small $\varepsilon$.
\par

Similarly,  for every
$\varphi\in C_0^\infty([0,+\infty)\times D)$ and for sufficiently small
$\varepsilon$,
\begin{equation}\label{Eq-4.1-10}
\int_0^t\langle \overline{\theta}_1^\varepsilon,
\varphi\rangle_{L^2(\partial D)}ds=O(\varepsilon),
\end{equation}
\begin{equation}\label{Eq-4.1-11}
\begin{array}{ll}
\int_0^t\langle \overline{\theta}_2^\varepsilon,
\varphi\rangle_{L^2(\partial D)}ds=&-\varepsilon \langle
\overline{\theta}_2^\varepsilon(t), \varphi(t)\rangle_{L^2(\partial
D)} +\varepsilon\int_0^t\langle \overline{\theta}_2^\varepsilon(s),
\varphi_s(s)\rangle_{L^2(\partial D)}ds\\
&-\int_0^t\langle \delta^\varepsilon, \varphi\rangle_{L^2(\partial
D)}ds+\int_0^t\langle v^\varepsilon, \varphi\rangle_{L^2(\partial
D)}ds,
\end{array}
\end{equation}
which implies from (\ref{Eq-p4.1-02}) that
\begin{equation}\label{Eq-4.1-12}
\begin{array}{ll}
\int_0^t\langle \overline{\theta}_2^\varepsilon,
\varphi\rangle_{L^2(\partial D)}ds=&O(\varepsilon)+O(\varepsilon)\\
&-\int_0^t\langle \delta^\varepsilon, \varphi\rangle_{L^2(\partial
D)}ds+\int_0^t\langle v^\varepsilon, \varphi\rangle_{L^2(\partial
D)}ds,
\end{array}
\end{equation}
and
\begin{equation}\label{Eq-4.1-13}
\begin{array}{ll}
\int_0^t\langle \overline{\theta}_3^{\varepsilon},
\varphi\rangle_{L^2(\partial D)}ds = &-\varepsilon \langle
\overline{\theta}_3^{\varepsilon}(t),
\varphi(t)\rangle_{L^2(\partial D)}ds+\varepsilon\int_0^t\langle
\overline{\theta}_3^{\varepsilon}(s),
\varphi_s(s)\rangle_{L^2(\partial D)}ds\\
&+\varepsilon^\alpha\int_0^t \langle dW_2(s),
\varphi\rangle_{L^2(\partial D)},
\end{array}
\end{equation}
which leads to
\begin{equation}\label{Eq-4.1-14}
\begin{array}{l}
\int_0^t\langle \overline{\theta}_3^{\varepsilon},
\varphi\rangle_{L^2(\partial D)}ds =
O(\varepsilon)+O(\varepsilon)+\varepsilon^\alpha\int_0^t \langle
dW_2(s), \varphi\rangle_{L^2(\partial D)}.
\end{array}
\end{equation}
Thus,  by the Gronwall
inequality, (\ref{Eq-p4.1-03}) and the condition $\alpha\in [1/2,
1)\cup (1, +\infty)$, we have
$\mathbb{E}\|\overline{\theta}_3^{\varepsilon}(t)\|_{L^2(\partial
D)}\leq TrQ_2$ for $t\in [0,+\infty)$.
\par
Therefore, substituting (\ref{Eq-4.1-2}), (\ref{Eq-4.1-5}),
(\ref{Eq-4.1-9}), (\ref{Eq-4.1-10}), (\ref{Eq-4.1-12}) and
(\ref{Eq-4.1-14}) into (\ref{Eq-4.1-1}), for every $\varphi\in
C_0^\infty([0,+\infty)\times D)$, we conclude that for sufficiently small
$\varepsilon$,
\begin{equation}\label{Eq-4.1-15}
\begin{array}{ll}
&\langle u^\varepsilon(t), \varphi(t) \rangle_{L^2(D)}+\langle
\delta^\varepsilon(t), \varphi(t) \rangle_{L^2(\partial
D)}\\
= &\langle u^\varepsilon(0), \varphi(0) \rangle_{L^2(D)}+\int_0^t
\langle u^\varepsilon(s), \varphi_s(s) \rangle_{L^2(D)}ds\\
&+\langle \delta^\varepsilon(0), \varphi(0) \rangle_{L^2(\partial
D)}+\int_0^t \langle \delta^\varepsilon(s), \varphi_s(s)
\rangle_{L^2(\partial D)}ds\\
 &+\int_0^t\langle \bigtriangleup u^{\varepsilon},
\varphi\rangle_{L^2(D)}ds- \int_0^t\langle u^{\varepsilon},
\varphi\rangle_{L^2(D)}ds+\int_0^t\langle \sin u^{\varepsilon},
\varphi\rangle_{L^2(D)}ds\\
&-\int_0^t\langle \delta^\varepsilon, \varphi\rangle_{L^2(\partial
D)}ds+\int_0^t\langle v^\varepsilon, \varphi\rangle_{L^2(\partial
D)}ds\\
& + \varepsilon^\alpha\int_0^t \langle dW_1(s),
\varphi\rangle_{L^2(D)}+\varepsilon^\alpha\int_0^t \langle
dW_2(s),\varphi\rangle_{L^2(\partial D)}+O(\varepsilon).
\end{array}
\end{equation}
Projecting (\ref{Eq-4.1-15}) onto $L^2(D)$ and $L^2(\partial D)$,
respectively, we   derive that for sufficiently small
$\varepsilon$, as $\alpha\in [1/2,1)$, the approximating equation for
$u^\varepsilon$  is
\begin{equation}\label{Eq-4.1-16}
\left\{
\begin{array}{ll}
\overline{u}^\varepsilon_t= \bigtriangleup
\overline{u}^{\varepsilon}-\overline{u}^{\varepsilon}+\sin
\overline{u}^{\varepsilon}+\varepsilon^\alpha\dot{\overline{W}_1},\quad
\overline{u}^\varepsilon(0)=u_0,&\quad \hbox{in}\quad D,\\
\overline{\delta}^\varepsilon_t=
-\overline{\delta}^{\varepsilon}+\overline{u}_t^{\varepsilon}+\varepsilon^\alpha\dot{\overline{W}_2},\quad
\overline{\delta}^\varepsilon(0)=\delta_0,&\quad \hbox{on}\quad
\partial D,\\
\overline{\delta}^\varepsilon_t=\frac{\overline{u}^\varepsilon}{\partial
{\bf n}},&\quad \hbox{on}\quad
\partial D.
\end{array}
\right.
\end{equation}
Then it follows from (\ref{Eq-4.1-15}) and (\ref{Eq-4.1-16}) that
the result under the condition $[1/2, 1)$ holds.
\par
It remains to consider   the case of $\alpha\in (1,
+\infty)$.
\par
It follows from (\ref{Eq-4.1-1}), (\ref{Eq-4.1-4}) and
(\ref{Eq-4.1-11}) that
\begin{equation}\label{Eq-4.1-17}
\begin{array}{ll}
&\langle u^\varepsilon(t), \varphi(t) \rangle_{L^2(D)} +\langle
\delta^\varepsilon(t), \varphi(t) \rangle_{L^2(\partial
D)}\\
= &\langle u^\varepsilon(0), \varphi(0) \rangle_{L^2(D)}+\int_0^t
\langle
u^\varepsilon(s), \varphi_s(s) \rangle_{L^2(D)}ds \\
&+\langle \delta^\varepsilon(0), \varphi(0) \rangle_{L^2(\partial
D)}+\int_0^t \langle \delta^\varepsilon(s), \varphi_s(s)
\rangle_{L^2(\partial D)}ds\\
&+ \int_0^t\langle \overline{v}_1^\varepsilon(s), \varphi(s)
\rangle_{L^2(D)}ds + \int_0^t\langle
\overline{v}_3^\varepsilon(s), \varphi(s) \rangle_{L^2(D)}ds\\
& + \int_0^t\langle \overline{\theta}_1^\varepsilon(s), \varphi(s)
\rangle_{L^2(\partial D)}ds  + \int_0^t\langle
\overline{\theta}_3^\varepsilon(s), \varphi(s) \rangle_{L^2(\partial
D)}ds\\
& -\varepsilon \langle \overline{v}_2^{\varepsilon}(t),
\varphi(t)\rangle_{L^2(D)}+ \varepsilon \int_0^t\langle
\overline{v}_2^{\varepsilon}(s),
\varphi_s(s)\rangle_{L^2(D)}ds\\
&+\int_0^t\langle \bigtriangleup u^{\varepsilon},
\varphi\rangle_{L^2(D)}ds- \int_0^t\langle u^{\varepsilon},
\varphi\rangle_{L^2(D)}ds+\int_0^t\langle \sin u^{\varepsilon},
\varphi\rangle_{L^2(D)}ds\\
&-\varepsilon \langle \overline{\theta}_2^\varepsilon(t),
\varphi(t)\rangle_{L^2(\partial D)} +\varepsilon\int_0^t\langle
\overline{\theta}_2^\varepsilon(s),
\varphi_s(s)\rangle_{L^2(\partial D)}ds\\
&-\int_0^t\langle \delta^\varepsilon, \varphi\rangle_{L^2(\partial
D)}ds+\int_0^t\langle v^\varepsilon, \varphi\rangle_{L^2(\partial
D)}ds.
\end{array}
\end{equation}
From (\ref{Eq-4-1}) and (\ref{Eq-4-2}), we have
\begin{equation}\label{Eq-4.1-18}
\begin{array}{l}
\overline{v}_2^{\varepsilon}=v^\varepsilon -\overline{v}_1^{\varepsilon}-\overline{v}_3^{\varepsilon},\\
\overline{\theta}_2^{\varepsilon}=\theta^\varepsilon
-\overline{\theta}_1^{\varepsilon}-\overline{\theta}_3^{\varepsilon}.
\end{array}
\end{equation}
Then, from (\ref{Eq-4.1-17}) and \ref{Eq-4.1-18}), we infer that
\begin{equation}\label{Eq-4.1-19}
\begin{array}{ll}
&\langle u^\varepsilon(t), \varphi(t) \rangle_{L^2(D)} -\langle
u^\varepsilon(0), \varphi(0) \rangle_{L^2(D)}-\int_0^t \langle
u^\varepsilon(s), \varphi_s(s) \rangle_{L^2(D)}ds \\
&+\langle \delta^\varepsilon(t), \varphi(t) \rangle_{L^2(\partial
D)}-\langle \delta^\varepsilon(0), \varphi(0) \rangle_{L^2(\partial
D)}-\int_0^t \langle \delta^\varepsilon(s), \varphi_s(s)
\rangle_{L^2(\partial D)}ds\\
&-\int_0^t\langle \bigtriangleup u^{\varepsilon},
\varphi\rangle_{L^2(D)}ds+\int_0^t\langle u^{\varepsilon},
\varphi\rangle_{L^2(D)}ds-\int_0^t\langle \sin u^{\varepsilon},
\varphi\rangle_{L^2(D)}ds\\
&+\varepsilon \langle v^\varepsilon(t), \varphi(t) \rangle_{L^2(D)}
-\varepsilon\int_0^t \langle v^\varepsilon(s), \varphi_s(s)
\rangle_{L^2(D)}ds\\
& +\int_0^t\langle \delta^\varepsilon, \varphi\rangle_{L^2(\partial
D)}ds-\int_0^t\langle v^\varepsilon, \varphi\rangle_{L^2(\partial
D)}ds\\
&+\varepsilon \langle \theta^\varepsilon(t), \varphi(t)
\rangle_{L^2(\partial D)}-\varepsilon \int_0^t \langle
\theta^\varepsilon(s), \varphi_s(s) \rangle_{L^2(\partial D)}\\
=&\varepsilon\langle \overline{v}_1^\varepsilon(t), \varphi(t)
\rangle_{L^2(D)}+\int_0^t\langle \overline{v}_1^\varepsilon(s),
\varphi(s) \rangle_{L^2(D)}ds-\varepsilon\int_0^t\langle
\overline{v}_1^\varepsilon(s), \varphi_s(s) \rangle_{L^2(D)}ds\\
&+\varepsilon\langle \overline{v}_3^\varepsilon(t), \varphi(t)
\rangle_{L^2(D)}+\int_0^t\langle \overline{v}_3^\varepsilon(s),
\varphi(s) \rangle_{L^2(D)}ds-\varepsilon\int_0^t\langle
\overline{v}_3^\varepsilon(s), \varphi_s(s) \rangle_{L^2(D)}ds\\
&+ \varepsilon \langle \overline{\theta}_1^\varepsilon(t),
\varphi(t)\rangle_{L^2(\partial D)}+\int_0^t\langle
\overline{\theta}_1^\varepsilon(s), \varphi(s)\rangle_{L^2(\partial
D)}ds-\varepsilon\int_0^t\langle \overline{\theta}_1^\varepsilon(s),
\varphi_s(s)\rangle_{L^2(\partial D)}ds\\
& +\varepsilon \langle \overline{\theta}_3^\varepsilon(t),
\varphi(t)\rangle_{L^2(\partial D)}+\int_0^t\langle
\overline{\theta}_3^\varepsilon(s), \varphi(s)\rangle_{L^2(\partial
D)}ds-\varepsilon\int_0^t\langle \overline{\theta}_3^\varepsilon(s),
\varphi_s(s)\rangle_{L^2(\partial D)}ds.
\end{array}
\end{equation}
\par
For (\ref{Eq-4-1.1}), we see that for
 $\varphi\in
C_0^\infty([0,+\infty)\times D)$,
\begin{equation}\label{Eq-4.1-20}
\int_0^t\langle \frac{d\overline{v}_1^{\varepsilon}}{ds},
\varphi\rangle_{L^2(D)}ds=-\frac{1}{\varepsilon}\int_0^t \langle
\overline{v}_1^{\varepsilon}, \varphi\rangle_{L^2(D)}ds.
\end{equation}
Noticing that $\overline{v}_1^{\varepsilon}(0)=v_0$
and that
$$
\int_0^t\langle \frac{d\overline{v}_1^{\varepsilon}}{ds},
\varphi\rangle_{L^2(D)}ds= \langle \overline{v}_1^{\varepsilon}(t),
\varphi(t)\rangle_{L^2(D)}-\langle \overline{v}_1^{\varepsilon}(0),
\varphi(0)\rangle_{L^2(D)}-\int_0^t\langle
\overline{v}_1^{\varepsilon}(s), \varphi_s(s)\rangle_{L^2(D)}ds,
$$
we deduce from (\ref{Eq-4.1-20}) that
\begin{equation}\label{Eq-4.1-21}
\begin{array}{ll}
& \varepsilon\langle \overline{v}_1^\varepsilon(t), \varphi(t)
\rangle_{L^2(D)}+\int_0^t\langle \overline{v}_1^\varepsilon(s),
\varphi(s) \rangle_{L^2(D)}ds-\varepsilon\int_0^t\langle
\overline{v}_1^\varepsilon(s), \varphi_s(s) \rangle_{L^2(D)}ds\\
= &\varepsilon \langle v_0, \varphi(0)\rangle_{L^2(D)}.
\end{array}
\end{equation}
Also, from (\ref{Eq-4.1-7}), (\ref{Eq-4.1-8}) and $\alpha\in
(1,+\infty)$,  we
have that for sufficiently small $\varepsilon$,
\begin{equation}\label{Eq-4.1-22}
\begin{array}{ll}
&\varepsilon\langle \overline{v}_3^\varepsilon(t), \varphi(t)
\rangle_{L^2(D)}+\int_0^t\langle \overline{v}_3^\varepsilon(s),
\varphi(s) \rangle_{L^2(D)}ds-\varepsilon\int_0^t\langle
\overline{v}_3^\varepsilon(s), \varphi_s(s) \rangle_{L^2(D)}ds\\
=&\varepsilon^\alpha\int_0^t \langle \varphi, dW_1(s)\rangle_{L^2(D)}\\
=&O(\varepsilon^\alpha).
\end{array}
\end{equation}
\par
Similarly, we  derive that
\begin{equation}\label{Eq-4.1-23}
\begin{array}{ll}
& \varepsilon\langle \overline{\theta}_1^\varepsilon(t), \varphi(t)
\rangle_{L^2(\partial D)}+\int_0^t\langle
\overline{\theta}_1^\varepsilon(s), \varphi(s) \rangle_{L^2(\partial
D)}ds-\varepsilon\int_0^t\langle
\overline{\theta}_1^\varepsilon(s), \varphi_s(s) \rangle_{L^2(\partial D)}ds\\
= &\varepsilon \langle \theta_0, \varphi(0)\rangle_{L^2(\partial
D)},
\end{array}
\end{equation}
and
\begin{equation}\label{Eq-4.1-24}
\begin{array}{ll}
&\varepsilon\langle \overline{\theta}_3^\varepsilon(t), \varphi(t)
\rangle_{L^2(\partial D)}+\int_0^t\langle
\overline{\theta}_3^\varepsilon(s), \varphi(s) \rangle_{L^2(\partial
D)}ds-\varepsilon\int_0^t\langle
\overline{\theta}_3^\varepsilon(s), \varphi_s(s) \rangle_{L^2(\partial D)}ds\\
=&O(\varepsilon^\alpha),
\end{array}
\end{equation}
for sufficiently small $\varepsilon$.
\par
Substituting (\ref{Eq-4.1-21})-(\ref{Eq-4.1-24}) into
(\ref{Eq-4.1-19}), for every $\varphi\in C_0^\infty([0,+\infty)\times
D)$, we have
that for sufficiently small $\varepsilon$,
\begin{equation}\label{Eq-4.1-25}
\begin{array}{ll}
&\langle u^\varepsilon(t), \varphi(t) \rangle_{L^2(D)} -\langle
u^\varepsilon(0), \varphi(0) \rangle_{L^2(D)}-\int_0^t \langle
u^\varepsilon(s), \varphi_s(s) \rangle_{L^2(D)}ds \\
&+\langle \delta^\varepsilon(t), \varphi(t) \rangle_{L^2(\partial
D)}-\langle \delta^\varepsilon(0), \varphi(0) \rangle_{L^2(\partial
D)}-\int_0^t \langle \delta^\varepsilon(s), \varphi_s(s)
\rangle_{L^2(\partial D)}ds\\
&-\int_0^t\langle \bigtriangleup u^{\varepsilon},
\varphi\rangle_{L^2(D)}ds+\int_0^t\langle u^{\varepsilon},
\varphi\rangle_{L^2(D)}ds-\int_0^t\langle \sin u^{\varepsilon},
\varphi\rangle_{L^2(D)}ds\\
&+\varepsilon \langle v^\varepsilon(t), \varphi(t) \rangle_{L^2(D)}
-\varepsilon\int_0^t \langle v^\varepsilon(s), \varphi_s(s)
\rangle_{L^2(D)}ds\\
& +\int_0^t\langle \delta^\varepsilon, \varphi\rangle_{L^2(\partial
D)}ds-\int_0^t\langle v^\varepsilon, \varphi\rangle_{L^2(\partial
D)}ds\\
&+\varepsilon \langle \theta^\varepsilon(t), \varphi(t)
\rangle_{L^2(\partial D)}-\varepsilon \int_0^t \langle
\theta^\varepsilon(s), \varphi_s(s) \rangle_{L^2(\partial D)}\\
=&\varepsilon \langle v_0, \varphi(0)\rangle_{L^2(D)}+ \varepsilon
\langle \theta_0, \varphi(0)\rangle_{L^2(\partial
D)}+O(\varepsilon^\alpha).
\end{array}
\end{equation}
Projecting (\ref{Eq-4.1-25}) onto $L^2(D)$ and $L^2(\partial D)$,
respectively, we   obtain that for sufficiently small
$\varepsilon$ and for $\alpha\in (1, +\infty)$, the approximating
equation of $u^\varepsilon$  is
\begin{equation}\label{Eq-4.1-26}
\left\{
\begin{array}{ll}
\varepsilon
\overline{u}^\varepsilon_{tt}+\overline{u}^\varepsilon_t-\bigtriangleup
\overline{u}^\varepsilon+\overline{u}^\varepsilon-\sin
\overline{u}^\varepsilon=0, &\quad in\; D,\\
\varepsilon\overline{\delta}^\varepsilon_{tt}+\overline{\delta}^\varepsilon_t+\overline{\delta}^\varepsilon
=-\overline{u}^\varepsilon_t, &\quad on\;
\partial D ,\\
\overline{\delta}^\varepsilon_t=\frac{\partial
\overline{u}^\varepsilon}{\partial {\bf n}}, &\quad on\;
\partial D,\\
\overline{u}^\varepsilon(0)=u_0, \overline{u}^\varepsilon_t(0)=u_1,
\overline{\delta}^\varepsilon(0)=\delta_0,
\overline{\delta}^\varepsilon_t(0)=\delta_1.
\end{array}
\right.
\end{equation}
\par
Therefore, it follows from (\ref{Eq-4.1-25}) and (\ref{Eq-4.1-26})
that the result under the condition $(1, +\infty)$ holds. This
completes the proof of Theorem 4.1.\hfill$\blacksquare$
\par

\end{document}